\documentclass{article}

\usepackage{geometry}
\usepackage{fancyhdr}

\usepackage{amsmath, amssymb}

\usepackage{graphicx}
\usepackage{color}
\definecolor{bg}{rgb}{0.89, 0.95, 0.71} % background e3f3b4
\definecolor{ac}{rgb}{0.50, 0.18, 0.41} % accent c491b6
\definecolor{rc}{rgb}{0.77, 0.57, 0.71} % rule 7f2f69
\usepackage{listings}
\lstdefinestyle{latex}{language=TeX,
                       backgroundcolor=\color{bg},
                       basicstyle=\small\ttfamily,
                       frame=leftline,
                       xleftmargin=1.4em, 
                       framexleftmargin=.8em}
\lstdefinestyle{cmdline}{
                         }
\usepackage{url}
\usepackage[pdftex, 
            bookmarks, 
            colorlinks=true,
            linkcolor=black, 
            plainpages = false,
            pdfpagemode = UseNone,
            pdfstartview = FitH, 
            citecolor = ac, urlcolor = ac, filecolor = ac]{hyperref}

\usepackage[para]{footmisc}
\makeatletter
\long\def\@makefntext#1{\leavevmode
\@makefnmark\nobreak
\hskip.05em\relax#1%
}
\makeatother
\makeatletter % from e$I_e$a pease at csli.stanford.edu
\def\section{\@startsection {section}{1}{\z@}{2.5ex plus .6ex minus 
    .2ex}{1.0ex plus .15ex}{\hspace*{-3em}\Large\bf\color{ac}}}
\def\subsection{\@startsection{subsection}{2}{\z@}{1.5ex plus .3ex minus 
   .1ex}{.2ex plus .1ex}{\hspace*{-3em}\bf\large\color{ac}}}
\makeatother
\newlength{\rulelength}
\title{Getting something out of \LaTeX{}}
\author{Jim Hef{}feron}
\usepackage{multirow}
\usepackage{wrapfig}
\usepackage{subfig}

\RequirePackage{lineno}
\modulolinenumbers[1]
\expandafter\def\expandafter\normalsize\expandafter{%
    \normalsize
    \setlength\abovedisplayskip{10pt}
    \setlength\belowdisplayskip{10pt}
    \setlength\abovedisplayshortskip{5pt}
    \setlength\belowdisplayshortskip{10pt}
\setlength{\jot}{10pt}    
}

\begin{document}

%\thispagestyle{empty}
% \maketitle
%\setlength{\unitlength}{1in}
%\begin{picture}(0,0)(0,0)
% \put(-0.5,-.84){\mbox{\color{bg}\rule{2.32in}{1.15in}}}
%\end{picture}
%\makeatletter

\begin{center}

{\large\bf Convergence Analysis for Computation of Coupled Advection-Diffusion-Reaction Problems}\\

\vspace*{0.4cm}

W. B. Dong$^1$, H. S. Tang$^{1,}$\footnote{Correspondence: Hasnong Tang, htang@ccny.cuny.edu}, Y. J. Liu$^2$

$^1$Civil Engineering Department, City College of New York, CUNY, New York, USA 

$^2$School of Mathematics, Georgia Institute of Technology, Atlanta, USA

\end{center}

\vspace*{0.2cm}

\begin{center}
\parbox{0.9\hsize}
{\small 
\noindent ABSTRACT 
{\sf \     
A study is presented on the convergence of the computation of coupled advection-diffusion-reaction equations. In the computation, the equations with different coefficients and even types are assigned in two subdomains, and Schwarz iteration is made between the equations when marching from a time level to the next one. The analysis starts with the linear systems resulting from the full discretization of the equations by explicit schemes. Conditions for convergence are derived, and its speedup and the effects of difference in the equations are discussed. Then, it proceeds to an implicit scheme, and a recursive expression for convergence speed is derived. An optimal interface condition for the Schwarz iteration is obtained, and it leads to ``perfect convergence", that is, convergence within two times of iteration. Furthermore, the methods and analyses are extended to the coupling of the viscous Burgers equations. Numerical experiments indicate that the conclusions, such as the ``perfect convergence, " drawn in the linear situations may remain in the Burgers equations' computation.  

\vspace{0.1cm} 

}
}
\end{center}

\vspace{0.5cm} 

\noindent {\bf 1. Introduction} 

\vspace{0.2cm} 

Advection, diffusion, and reaction are fundamental physical phenomena, which and whose interactions take place in real-world problems. These phenomena and interactions are described by  advection-diffusion-reaction equations and their coupling, which have been solved numerically to understand such phenomena underlying various physical processes, such as migration of pollutants in a porous medium, chemical reaction in a nuclear reactor, and propagation of tsunamis in the ocean \cite {Atis2012, Hamilton2016, Qu2019}. Additionally, more complicated partial differential equations (PDEs) and their coupling, such as the Navier-Stokes equations and their hydrostatic versions \cite {Tang2014, Blayo2016}, are adopted to simulate realistic multiscale/multiphysics problems. Since advection-diffusion-reaction equations exhibit behaviors of parabolic, hyperbolic, and elliptic equations, they are commonly used as models to study these complicated PDEs and their coupling \cite {Gastaldi1989, Sakamoto2015, Main2018, Manzanero2020}. As a result, a study on the computation of coupled advection-diffusion-reaction equations promotes simulations of actual physical problems and the development of numerical methods for PDEs. 

Domain decomposition (DD) has emerged as an indispensable avenue to the scientific computation of various problems, and its development in methods and analysis for coupled equations of advection, diffusion, and reaction can trace back over 30 years ago, e.g., \cite {Gastaldi1989, Gander2008}. A natural approach of computation is to compute these equations in subdomains in conjunction with Schwarz iteration in between when marching from a time level to the next level \cite {Meurant1991, Hoang2017, Canuto2019}. This approach is widely used in practical problems such as simulation of fluid flows \cite {Dolean:2002do,Tang2003}. Another popular approach is to adopt Schwarz waveform relaxation, by which at each time of Schwarz iteration the equations in subdomains are solved at all time levels \cite {Gander1998, Martin2004, Califano2018}. Investigations have covered diffusion equations \cite {Meurant1991, Gander1998, Linel2013}, advection equations \cite {Bamberger1997, Dolean:2002do,Dolean2008}, advection-diffusion equations \cite {Giladi2002,Martin2004,Gander2007a}, heterogeneous advection-diffusion equations \cite {Gander2007b,Halpern2012}, and other closely-related equations (e.g., the Darcy law, a fractional equation, and the Schrödinger equation) \cite {Hoang2017, Califano2018, Antoine2019}, all with or without reaction. The research topics include, algorithm study \cite {Daoud2002, Zhu2009, Canuto2019}, convergence analysis \cite {Dolean2008,Eisenmann2018, Califano2018, Canuto2019}, interface conditions in heterogeneous problems (e.g., between an advection-diffusion equation and an advection equation) \cite {Gastaldi1989, Gander2009, Hoang2017}, and optimal transmission conditions \cite {Gander2007a, Martin2004}. As a more general form of  advection-reaction equations, hyperbolic systems of conservation laws have also been investigated intensively \cite {Starius1980, Berger1987}, especially in the context of fluid flow simulation \cite {Tang2020}. In the computation of the systems, explicit discretization is commonly used, and, therefore, no iteration between subdomains is needed. In such a situation, interface algorithms and stability issues have been a main research topic \cite {Berger1987, Part1994, Tang1999}. 

Progress has been made in understanding DD computation of coupled advection-diffusion-reaction equations.  
%It is concluded that, when waveform relaxation is in use for an ODE, which may result from the spatial discretization of a PDE, the convergence rate is linear when an unbounded time interval is considered, and it may become super-linear when the interval is bounded \cite {Gander1998}. 
%Furthermore, 
%It is the convergence rate gets slower as the grid spacing gets finer \cite {Gastaldi1989, Gander1998}. 
It is known that the convergence rate increases with the size of the overlapping region,  decreases with diffusion coefficients, and gets slower as the grid spacing gets finer \cite {Gander2007a}. 
%It is agreed that convergence is slow when a regular Dirichlet transmission condition is adopted, and its speed up is possible by adopting an optimal condition \cite {Martin2004, Califano2018}. 
The convergence rate can be accelerated by methods such as precondition in iteration matrices and optimization for interface conditions, and such speedup can be over an order of magnitude \cite {Martin2004, Gander2008}. For systems of conservation laws, it is proven that, when conservative interface algorithms are adopted at grid interfaces, numerical solutions converge to weak solutions of the systems if they converge \cite {Berger1987}. Then, it is shown that, under certain conditions that can be verified, even when non-conservative interface algorithms are used, still convergent solutions will become weak solutions, and the conservation errors are bounded \cite {Tang1999}. Later, it is further proved that the conservation error tends to zero at the speed of $\cal O$ ($\Delta t$), regardless of the smoothness of solutions at the interfaces \cite {Tang2007}. In the computation of the systems, discretion is necessary to avoid numerical oscillations, nonphysical solutions, and multiple solutions \cite {Part1994, Wu1996, Tang1996}.

Despite the substantial study, more topics deserve investigation. Heterogeneous DD problems of advection-diffusion-reaction equations are such a topic. In the past, research has been focused on homogeneous problems, e.g., identical equations in subdomains, and investigations on heterogeneous problems are relatively sparse, although relevant efforts have been made previously \cite {Gastaldi1989, Gander2009}. Additionally, previous analysis investigations are primarily conducted at the continuous level and semi-discretizations levels \cite {Houzeaux2003,Eisenmann2018}. While studies at such levels are successful and have revealed many important patterns of the problems, their results may not be conveniently or directly applied to practical computation. For instance, often the analysis of interface conditions is made based on the Laplace or Fourier transform, and the corresponding obtained optimal parameters are expressed in a spectral space and thus need to be approximated for actual computation \cite {Dolean:2002do, Gander2004}. Moreover, most studies are about linear problems in the past, with scarce discussion on nonlinear equations. 
%Therefore, coupling of nonlinear advection-diffusion-reaction problems is an important topic.  

This paper presents a systematic analysis of the computation of coupled advection-diffusion-reaction equations, and it intends to fill the gaps listed above and also deals with situations beyond those in existing work. Notably, such equations, with different coefficients in subdomains in general, are fully discretized using either explicit or implicit schemes. The analysis starts with explicit schemes for the computation of the resulting algebraic systems. Conditions for the convergence are derived, and its speedup and behaviors are examined, together with a discussion on the effects of difference in PDEs in subdomains. Then, the analysis proceeds to an implicit scheme and derivation of its convergence speed, followed by the derivation of its optimal interface algorithm. Finally, the study is extended to the coupling of the viscous Burgers equations, and analysis based on linearization is made on their computation.  

\vspace{0.3cm} 

\noindent {\bf 2. The problem of study}

\vspace{0.2cm} 

This paper considers computation of the following initial value problem of coupled advection-diffusion-reaction equations in two subdomains: 
\begin {equation}
\label{ibvp_adr}
\left\{\begin{array}{ll}
{u_1}_t+f^1_x(u_1)=b_1{u_1}_{xx}-c_1u_1, \   & x<x_2 \\
u_1=u_2, \ & x=x_2  \\
u_1=g(x), \  & t=0  
\end{array}\right. \quad  
\left\{\begin{array}{ll}
{u_2}_t+f^2_x(u_2)=b_2{u_2}_{xx}-c_2u_2,  \  &x>x_1 \\
u_2=u_1, \ & x=x_1\\
u_2=g(x), \   & t=0  \\
\end{array}\right. 
\end {equation} 
\noindent in which $t$ is the time, and $x$ is the space coordinate. Subdomain 1 is on the left, $x\le x_2$, and subdomain 2 is on the right, $x\ge x_1$. The two subdomains overlap with each other, i.e., $x_1<x_2$, and $x=x_1$ and $x_2$ are the interface of subdomain 2 and 1, respectively. In this study, $f^k(u_k)=a_ku_k, {u_k}^2/2$ ($k=1,2$) will be considered. $a_k$, $b_k$, and $c_k$ are constant coefficients, and their values may be different in the two subdomains. As a result, problem (\ref {ibvp_adr}) includes a number of coupling scenarios, such as coupling between two advection-diffusion-reaction equations with different coefficient, between an advection-diffusion equation ($a_1,b_1,c_1\ne 0$) and an advection equation ($b_2,c_2=0$), and an advection equation ($b_1,c_1=0$) and a diffusion-reaction equation ($a_2=0$). As seen in (\ref {ibvp_adr}), Dirichlet conditions are used the interfaces. Moreover, boundary conditions at the left and right ends of subdomain 1 and 2, respectively, may be added. Without loss of generality, it is assume that $g(x)$ has a finite support. 

\begin{figure}[ht]
\centering
\includegraphics[scale=0.6]{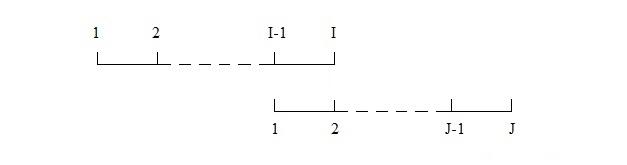}
\caption{\small Computational grids. The left and right grids cover subdomain 1 and 2, respectively.}
\label{dd_grids}
\end{figure}

Consider computation of problem (\ref {ibvp_adr}) on two grids as shown in Fig.\ref {dd_grids}, with the left and right grid covering subdomain 1 and 2, respectively. Node $I$ of the left grid and node $1$ of the right grid are the  interfaces of the two grids, and they correspond to the interfaces of the two subdomains. On the two grids, the advection-diffusion-reaction equations are discretized with backward difference and central difference in time and space, respectively, and this leads to an implicit scheme:
\begin{equation}
\label{adr_dis_1}
\frac{{u_k}_i^{n+1}-{u_k}_i^n}{\Delta t} + {f_u^k}_i^{n+1}\frac{{u_k}_{i+1}^{n+1}-{u_k}_{i-1}^{n+1}}{2\Delta x}=b_k\frac{{u_k}_{i+1}^{n+1}-2{u_k}_i^{n+1}+{u_k}_{i-1}^{n+1}}{\Delta x^2}-c_k{u_k}_i^{n+1}
\end{equation} 

\noindent where $k=1,2$, depicting subdomain 1 and 2, respectively, $\Delta t$ and $\Delta x$ are the time step and the grid spacing, respectively, superscript $n$ and subscript $i$ indicate time level and grid node, respectively, and $f_u^k=\partial f^k(u_k)/\partial u_k$.  

This paper aims to analyze convergence of the computation of discretization (\ref {adr_dis_1}). The computation can be carried out by treating ${u_k}_i^{n+1}$ as ${u_k}_i^m$ and/or ${u_k}_i^{m+1}$, which leads to either explicit or implicit schemes (in terms of ${u_k}_i^{m+1}$) in the two subdomains, together with Schwarz iteration in between. Particularly, when marching from time level $n$ to $n+1$ via an explicit scheme, the problem  will be solved  as ($m=0,1,2, ...$): 
%{\bf double check the following interface condition?}
\begin{equation}
\label{Schwarz_iteration_1}
\left\{\begin{array}{ll} 
\mathcal {F}(...,{u_1}_i^m,...,{u_1}_i^{m-1},...,{u_1}_i^n,...)=0, \
& \ {u_1}_I^m={u_2}_1^m  \\ 
\mathcal {F}(...,{u_2}_i^m,...,{u_2}_i^{m-1},...,{u_2}_i^n,...)=0, \ 
& \ {u_2}_1^m={u_1}_{I-1}^m
\end{array}\right .
\end{equation}
\noindent Here, $\mathcal {F}$ represents an iterative scheme to compute discretization (\ref {adr_dis_1}) in subdomains. As $m\rightarrow \infty $, if convergent, the iterated solution converges to the solution at time level $n+1$: ${u_1}_i^m \rightarrow {u_1}_i^{n+1}$, ${u_2}_i^m\rightarrow {u_2}_i^{n+1}$. In above, the two grids exchange solutions at a same level of iteration, i.e. $m$, at their interfaces, and this is possible when an explicit iterative scheme is adopted in suddomains. When an implicit scheme is used, the exchange will be made between different levels of iteration: ${u_1}_I^m={u_2}_1^{m-1}$, and ${u_2}_1^m={u_1}_2^{m-1}$. Such interface conditions, with solution exchange either at a same or at a different level, are Dirichlet conditions, and they are referred to as classic interface conditions hereafter. 

Note that, instead of (\ref {adr_dis_1}), the equations may be discretized using an explicit scheme, but this study focuses on an implicit scheme given its advantages such as better numerical stability and wide use in practical computation. Moreover, it is readily checked that the discretization accuracy of (\ref {adr_dis_1}) is ${\cal O} (\Delta t+\Delta x^2)$. If needed, ${\cal O} (\Delta t^2+\Delta x^2)$ can be achieved, for instance, by discretizing the time derivative with a second-order accurate, three-time-level  backward difference (e.g., \cite {Tang2003}), or, discretizing the whole equation with the Crank-Nicholson scheme (e.g., \cite {Tang2007}). However, the analysis of (\ref {adr_dis_1}) in this paper remains the same or similar when these types of discretization are considered. At last, instead of that by waveform relaxation, this work considers the computation by Schwarz iteration only within two adjacent time levels $n$ and $n+1$, which is commonly adopted, especially in practical problem \cite {Tang2003, Canuto2019}. 

\vspace{0.3cm} 

\noindent {\bf 3. Explicit scheme}

\vspace{0.2cm} 

\noindent {\bf 3.1 Preliminary}

\vspace{0.2cm}

Consider a linear system 
\begin {equation}
\label{linear_sys_1}
AU=b
\end {equation} 
\noindent in which $A=(a_{ij})$, $U=(u_i)$, and $b$=$(b_i)$, with $i,j=1,2, ..., N$. When the system is solved by the Jacobi iteration, one has ($m=1, 2, ... $)
\begin {equation}
\label{linear_sys_2}
U^m=(I-D^{-1}A)U^{m-1}+D^{-1}b
\end {equation} 
\noindent in which $I$ is the identity matrix, and $D$ is a diagonal matrix with the diagonal elements of $A$ as its diagonal elements. A theorem shown by  Saad in \cite {Saad2003} shows that (\ref {linear_sys_2}) converges as long as 
\begin{equation}
\label{converge_cond2}
\sum_{i=1,\neq j}^N |a_{ij}|/|a_{ii}|
\left\{\begin{array}{ll}
\le 1, \ j=1,..., N \\
<1,  \ at \ least \ one \ j\ 
\end{array}\right . 
\end{equation}

\noindent and $A$ is irreducible. Instead of above, the following condition is frequently used as a sufficient condition for convergence of (\ref {linear_sys_2}) in practical computation, e.g., heat and mass transfer \cite{patankar:1980,versteeg:2007}: 
\begin{equation}
\label{Scarborough_cond}
\sum_{j=1,\neq i}^N |a_{ij}|/|a_{ii}|
\left\{\begin{array}{ll}
\le 1, \ i=1,..., N \\
<1,  \ at \ least \ one \ i\ 
\end{array}\right . 
\end{equation}
\noindent which is proposed by Scarborough \cite {James:1930,James:1958}, and is referred to as the Scarborough criterion \cite{patankar:1980,versteeg:2007}. Note that, the summation in condition (\ref {converge_cond2}) is over rows, while that in the Scarborough criterion is over columns. 

{\it Remark 3.1} In the work of Scarborough \cite {James:1958}, criterion (\ref {Scarborough_cond}) is proposed but not actually proved, and in fact condition (\ref {converge_cond2}) is proved instead, all without "=". The theorem shown by Saad in \cite {Saad2003} (Theorem 4.9) actually leads to a proof for the Scarborough criterion. Actually, by the theorem, Scarborough criterion (\ref {Scarborough_cond}) is a convergence condition for iteration $U^m=(I-D^{-1}A)^TU^{m-1}$. As a result, the spectral radius of $(I-D^{-1}A)^T$ is less than $1$, and thus the spectral radius of $(I-D^{-1}A)$ must be less than $1$ also since their radius are the same. As a result, iteration (\ref{linear_sys_2}) converges under the Scarborough criterion.    

The following is a trivial case for a more general conclusion. Given an iteration 
\begin {equation}
\label{linear_sys_3}
U^m=CU^{m-1}+d
\end {equation}  

\noindent where $C=(c_{ij})$, and $d=(d_i)$, with $i,j=1,2, ..., N$. Let 
\begin {equation}
\label{K}
K=\max_i\bigg\{ \sum_{j=1}^{L} |c_{ij}| \bigg\}=||C||_\infty
\end {equation} 

\noindent In view that $||C||<1$ leads to the convergence of the iteration, $K<1$, or, the following 
\begin{equation}
\label{converge_cond3}
\sum_{j=1}^{N} |c_{ij}| <1
\end{equation}
for every $i$, is a condition for the convergence.

\vspace{0.3cm} 

\noindent {\bf 3.2 Condition for convergence}

\vspace{0.2cm} 

Consider computation of (\ref {Schwarz_iteration_1}) by iterative, explicit schemes in the two subdomains. Two common explicit schemes will be considered, and the first one reads as ($k=1,2$)
\begin{equation}
\label{adr_dis_2}
\dfrac{{u_k}_i^m-{u_k}_i^n}{\Delta t} + a_k\dfrac{{u_k}_{i+1}^{m-1}-{u_k}_{i-1}^{m-1}}{2\Delta x}=
b_k\dfrac{{u_k}_{i+1}^{m-1}-2{u_k}_i^m+{u_k}_{i-1}^{m-1}}{\Delta x^2}-c_k{u_k}_i^m,  
\end{equation} 
\noindent which is iterative and actually the Jacobi iteration \cite {Epperson:2013}. This scheme is simple and slow in convergence, but it permits parallel computation and is in frequent use for practical problems \cite{patankar:1980,versteeg:2007}. With this scheme, Schwarz iteration (\ref {Schwarz_iteration_1}) reads as   
\begin{eqnarray}
%\begin {equation}
\label{Schwarz_iteration_2}
\left\{\begin{array}{lll}
(1+2\epsilon  _1+\gamma _1){u_1}_i^m=(\eta _1+\epsilon _1){u_1}_{i-1}^{m-1}+(\epsilon  _1-\eta  _1){u_1}_{i+1}^{m-1}+{u_1}_i^n, & i\le I-1;  \  
&{u_1}_I^m={u_2}_2^m \\
(1+2\epsilon  _2+\gamma _2){u_2}_i^m=(\eta  _2+\epsilon _2){u_2}_{i-1}^{m-1}+(\epsilon  _2-\eta  _2){u_2}_{i+1}^{m-1}+{u_2}_i^n, & i\ge 2; \ &{u_2}_1^m={u_1}_{I-1}^m 
\end{array}\right.  
%\end {equation} 
\end{eqnarray}

\noindent in which $\eta_k =a_k\Delta t/(2 \Delta x), \epsilon_k =b_k\Delta t/\Delta x^2$, and $\gamma _k=\Delta tc_k $ ($k=1,2$). By arranging terms, above iteration can be expressed as (\ref {linear_sys_2}), in which 
\begin {equation*}
\label{matrix_1a}
\begin{array}{lll}
A=
   \begin{bmatrix} 
   \begin{array}{cccc|cccc}
 \ddots  & \ddots & \ddots &        &        &          &          &\\
         & r_1    & s_1    & t_1    &        &          &          &\\
         &        & r_1    & s_1    & t_1    &          &          &\\
  \hline
         &        &        & r_2    & s_2    & t_2     &           &\\
         &        &        &        & r_2    & s_2     &   t_2     &\\    
         &        &        &        &        & \ddots  &   \ddots  & \ddots 
   \end{array}
   \end{bmatrix}, \ \ 
U= 
   \begin{bmatrix} 
   \begin{array}{c}
   \vdots\\
   {u_1}_{I-2} \\
   {u_1}_{I-1} \\
  \hline
   {u_2}_2 \\  
   {u_2}_3 \\
   \vdots
   \end{array}
   \end{bmatrix}, 
\  \ 
b=
   \begin{bmatrix} 
   \begin{array}{c}
   \vdots\\
   {u_1}_{I-2}^n \\
   {u_1}_{I-1}^n \\
  \hline 
   {u_2}_{2}^n \\
   {u_2}_{3}^n \\
   \vdots
   \end{array}
   \end{bmatrix}    
\end{array}  
\end {equation*}

\noindent in which $r_k=-\epsilon_k-\eta_k, \ s_k=1+2\epsilon_k+\gamma _k, \ t_k=\eta_k-\epsilon_k$. By Scarborough principle (\ref {Scarborough_cond}), the following is concluded.

{\textsc{proposition 3.1}} A sufficient condition for iteration (\ref {Schwarz_iteration_2}) to converge is
\begin {equation}
\label{convergence_cond_1}
\max\limits_k\bigg\{\frac {|\epsilon_k+\eta_k |+|\eta_k-\epsilon_k|}{|1+2\epsilon_k+\gamma _k|}\bigg\}<1
\end {equation}

Now, let us discuss if (\ref {convergence_cond_1}) is, or, how to make it satisfied. In view that $\eta_k =a_k\Delta t/(2 \Delta x)$, $\epsilon_k =b_k\Delta t/\Delta x^2$, and $\gamma _k=\Delta tc_k $, it is readily seen that, the LHS of (\ref {convergence_cond_1}) decreases with $\Delta t$, and, as $\Delta t$ becomes sufficiently small, (\ref {convergence_cond_1}) will be guaranteed. Additionally, it is seen that the LHS also decreases with $\gamma _k$. Note that the LHS is $K$, or, $||C||_\infty $, see (\ref {K}). In view that $\rho (C)\le ||C||_\infty$, a small value of $\Delta t$ and a large value of $\gamma _i$ lead to a smaller $K$ and thus are helpful to assure the convergence. 

The coefficients of the PDEs may be different in the two subdomains. For instance, when  $\epsilon _1, \gamma _1, \gamma _2=0$, the coupling becomes one between an advection equation and a advection-diffusion equation. In this scenario, the LHS becomes
\begin {equation*}
K=\mbox {max} \bigg\{ \bigg| \frac {2\eta _1}{1+2\epsilon_1}\bigg|,\frac {|\epsilon_2+\eta_2 |+|\eta_2-\epsilon_2|}{|1+2\epsilon_2|} \bigg\}
\end {equation*}

\noindent which will be less than 1 as long as $\Delta t$ is sufficiently small, as discussed above. Another scenario is $\eta _k=0$, that is, the problem becomes coupling between two diffusion equations. Because $\epsilon _k,\gamma _k>0$, one has
\begin {equation*}
K=\max\limits_k\bigg| \frac {2\epsilon _k}{1+2\epsilon_k+\gamma _k}\bigg| <1
\end {equation*}

\noindent or, (\ref {convergence_cond_1}) holds automatically. Discussion for more scenarios can be made. 
%In summary, condition (\ref {convergence_cond_1}) holds as time step becomes sufficiently small, or, automatically. 

Another iterative, explicit algorithm is the so-called artificial compressibility method, which is a main approach used in computation of incompressible flow problems \cite {Ge2005a, Tang2007}. By this approach, an artificial term is added, and the discretization reads as  
\begin{equation}
\label{adr_dis_3}
\frac{{u_k}_i^m-{u_k}_i^{m-1}}{\Delta \tau}+\frac{{u_k}_i^m-{u_k}_i^n}{\Delta t} + a_k\frac{{u_k}_{i+1}^{m-1}-{u_k}_{i-1}^{m-1}}{2\Delta x}=b_k\frac{{u_k}_{i+1}^{m-1}-2{u_k}_i^{m-1}+{u_k}_{i-1}^{m-1}}{\Delta x^2}-c_k{u_k}_i^{m-1}
\end{equation} 
\noindent in which $\Delta \tau$ is a pseudo-time step. At the convergence, the artificial term disappears. While the artificial compressibility method is adopted, the Schwarz iteration becomes 
\begin{eqnarray}
\label{Schwarz_iteration_3}
\left\{\begin{array}{lll}
(1+\kappa ){u_1}_i^m
=(\eta _1+\epsilon _1){u_1}_{i-1}^{m-1}
+\left (\kappa -2\epsilon  _1-\gamma _1\right ){u_1}_i^{m-1}
+(\epsilon  _1-\eta  _1){u_1}_{i+1}^{m-1}+{u_1}_i^n,   \\
\qquad i\le I-1; \quad {u_1}_I^m={u_2}_2^m \\
(1+\kappa ){u_2}_i^m
=(\eta _2+\epsilon _2){u_2}_{i-1}^{m-1}
+\left (\kappa -2\epsilon  _2-\gamma _2\right ){u_2}_i^{m-1}
+(\epsilon  _2-\eta  _2){u_2}_{i+1}^{m-1}+{u_2}_i^n,   \\
\qquad i\ge 2; \quad {u_2}_1^m={u_2}_{I-1}^m
\end{array}\right.  
\end{eqnarray}
\noindent where $\kappa=\Delta t/\Delta \tau$. Above can be formulated into form of iteration (\ref {linear_sys_3}). In this situation, $C$ is a tri-diagonal matrix, and, from the left to the right in a row, its non-zero elements are 
$   r_k=(\epsilon_k+\eta_k)/(1+\kappa ), \ 
    s_k=(\kappa -(2\epsilon_k+\gamma _k))(1+\kappa ), \ 
    t_k=(\epsilon_k-\eta_k)(1+\kappa )
$, together with   
$d=({u_1}_{2}^n/(1+\kappa),\cdots,{u_1}_{I-1}^n/(1+\kappa),{u_2}_{2}^n/(1+\kappa),\cdots,{u_2}_{J-1}^n/(1+\kappa))^T$. According to (\ref {converge_cond3}), we have the following conclusion. 

{\textsc{proposition 3.2}} A sufficient condition for (\ref {Schwarz_iteration_3}) to converge is 
\begin {equation}
\label{convergence_cond_2}
\max\limits_k\bigg\{
\bigg| \dfrac {\epsilon_k+\eta_k}{1+\kappa } \bigg|+\bigg| \dfrac{\kappa -(2\epsilon_k+\gamma _k)}{1+\kappa }\bigg|+\bigg|\dfrac {\epsilon_k-\eta_k}{1+\kappa }\bigg|\bigg \}
<1
\end {equation}

Now, let us also examine if (\ref {convergence_cond_2}) can be satisfied. Actually, (\ref {convergence_cond_2}) will hold as long as $\Delta \tau $ is sufficiently small. This is because in this situation, when $\epsilon _k\ge \eta _k $, its LHS becomes 
\begin {equation*}
K=\max\limits_k \bigg|\dfrac {\kappa -\gamma _k}{1+\kappa }\bigg|<1
\end {equation*}
\noindent When $\epsilon _k<\eta _k$, it is easy to check that 
\begin {equation*}
K=\max\limits_k\bigg|\dfrac {\kappa -(2(\epsilon _k-\eta _k)+\gamma _k)}{1+\kappa }\bigg|<1
\end {equation*}
\noindent Discussions can be made for more possible scenarios. 
%In summary, condition (\ref {convergence_cond_2}) holds as the pseudo time step is sufficiently small, or, automatically. 

\vspace{0.3cm} 

\noindent {\bf 3.3 Analysis of convergence speed}

\vspace{0.2cm} 

\underline {Convergence speed} Introducing relaxation into the interface condition in (\ref {Schwarz_iteration_2}), one has   
\begin{eqnarray*}
\label{relaxation_int}
\begin{array}{ll}
{u_1}_I^m={u_2}_2^{m-1}+\omega _1 ({u_2}_2^m-{u_2}_2^{m-1}) \\
{u_2}_1^m={u_1}_{I-1}^{m-1}+\omega _2 ({u_1}_{I-1}^m-{u_1}_{I-1}^{m-1}) 
\end{array}
\end{eqnarray*}
\noindent where $\omega _1, \omega _2=consts$, being the relaxation coefficients. With such an interface condition, the Jacobi iteration (\ref {linear_sys_2}) is formulated as (\ref {linear_sys_3}), in which 
\begin {equation*}
%\label{matrix_Jacobi_classic_relax}
\begin{array}{lll}
C=
   \begin{bmatrix} 
   \begin{array}{cccc|cccc}
 \ddots  & \ddots & \ddots &        &        &          &          &\\
         & r_1    & 0    & t_1    &        &          &          &\\
         & 0      & 0    & 0    & \omega _1 r_2    & (1-\omega _1)     & \omega_1 t_2    &\\
  \hline
         & \omega_2 r_1 &  (1-\omega_2) & \omega_2 t_1   & 0    & 0     &   0        &\\
         &        &        &        & r_2    & 0     &   t_2     &\\    
         &        &        &        &        & \ddots  &   \ddots  & \ddots 
   \end{array}
   \end{bmatrix}, \ \
U= 
   \begin{bmatrix} 
   \begin{array}{c}
   \vdots\\
   {u_1}_{I-1} \\
   {u_1}_I \\
  \hline
   {u_2}_1 \\  
   {u_2}_2 \\
   \vdots
   \end{array}
   \end{bmatrix} 
\  \ 
\end{array}  
\end {equation*}
\noindent and
$  r_k=(\epsilon_k+\eta_k)/(1+2\epsilon_k+\gamma_k), \ 
   t_k=(\epsilon_k-\eta_k)/(1+2\epsilon_k+\gamma_k)
$. According to condition (\ref {converge_cond3}), a sufficient condition for convergence is 
\begin{eqnarray*}
\label{Jacobi_classic_relax_convergence_cond}
\max\limits_{k,k'}\bigg\{|\omega _{k} |\left ( \Big| \dfrac {\epsilon_{k'}+\eta _{k'}}{1+2\epsilon_{k'}+\gamma _{k'}}\Big| 
               + \Big| \dfrac {\epsilon_{k'}-\eta _{k'} }{1+2\epsilon_{k'}+\gamma _{k'}}\Big| \right ) 
+|1-\omega _k |\bigg\}< 1 
\end{eqnarray*}
\noindent where $k,k'=1,2$, and $k\ne k'$. Note that the LHS becomes the summation of the other rows without $\omega _k, \omega _{k'}$ in $C$ when $\omega _k, \omega _{k'}=1$. As discussed above, it is anticipated that the convergence of the iteration may become faster as $K$ becomes smaller. It is readily checked that $\partial K/\partial \omega _k<-1$, $>-1$, and $>1$ when $\omega _k\le 0$, $0<\omega _k\le 1$, and $\omega _k>1$, respectively. It is expected that $K$ will have a minimum within $0\le \omega _k \le 1$, and, particularly, at either $\omega _k=0$, or, $=1$, which correspond to interface condition ${u_1}_I^m={u_2}_2^{m-1},\ {u_2}_1^m={u_1}_{I-1}^{m-1} $ and that in (\ref {Schwarz_iteration_2}), respectively.  

As a numerical experiment, three cases shown in Table \ref {linear_cases} are computed. The three cases represent typical situations of problem (\ref {ibvp_adr}). Particularly, Case 1, 2, and 3 involve coupling between two diffusion-dominant equations, between two advection-dominant equations, and between a diffusion-dominant and an advection-dominant equation, respectively. The result is given in Table \ref {diff_same_level}, which shows that spectral radius, $\rho $, for $\omega _k=1$ is lower than that for $\omega _k=0$, exhibiting a faster convergence.   
\begin{table}[ht]
  \centering
  \caption{\small Cases for numerical experiments on linear advection-diffusion-reaction equations. } 
  \label{linear_cases}
  \begin{tabular}{|c|c|c|c|c|c|c|}
\hline
       Case  & $a_1$ & $b_1$ & $c_1$ & $a_2$ & $b_2$ & $c_2$ \\ \hline
        1    & 0.01  & 0.5    & 0.0    & 0.01  & 0.5    & 0.0   \\ \hline       
        2    & 0.5   & 0.002  & 0.0    & 0.25  & 0.001  & 0.0   \\ \hline                  
        3    & 0.25  & 0.01   & 0.25   & -0.25 & 0.01   & 0.5   \\ \hline        
  \end{tabular}
\end{table}

\begin{table}[ht]
\centering
\caption{\small Spectral radii of the Jacobi iteration associated with the relaxation interface condition. (I,J)=(40,80).} 
\begin{tabular}{|c|c|c|c|}
\hline
Case          &  $\omega _1, \omega _2=0$   &  $\omega _1,\omega _2=1$ \\ \hline
1             & 0.967363 &   0.966594 \\ \hline
2             & 0.267145 &   0.197162 \\ \hline
3             & 0.431927 &   0.359324 \\ \hline
\end{tabular}
\label{diff_same_level}
\end{table}

\begin{figure}[ht] 
\centering
\subfloat[][]{\label{}\includegraphics[width=0.4\linewidth]{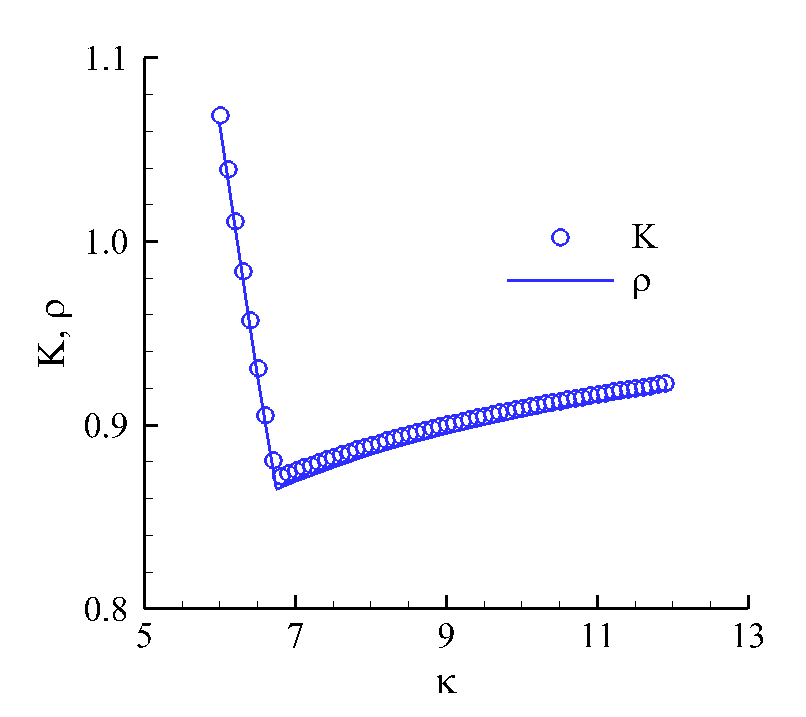}} 
\subfloat[][]{\label{}\includegraphics[width=0.4\linewidth]{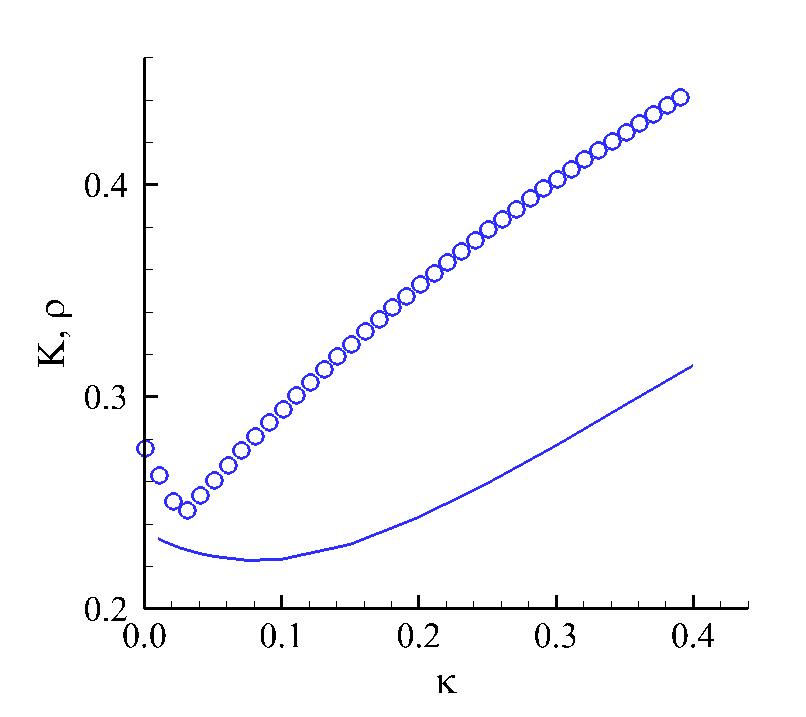}} \\
\subfloat[][]{\label{}\includegraphics[width=0.4\linewidth]{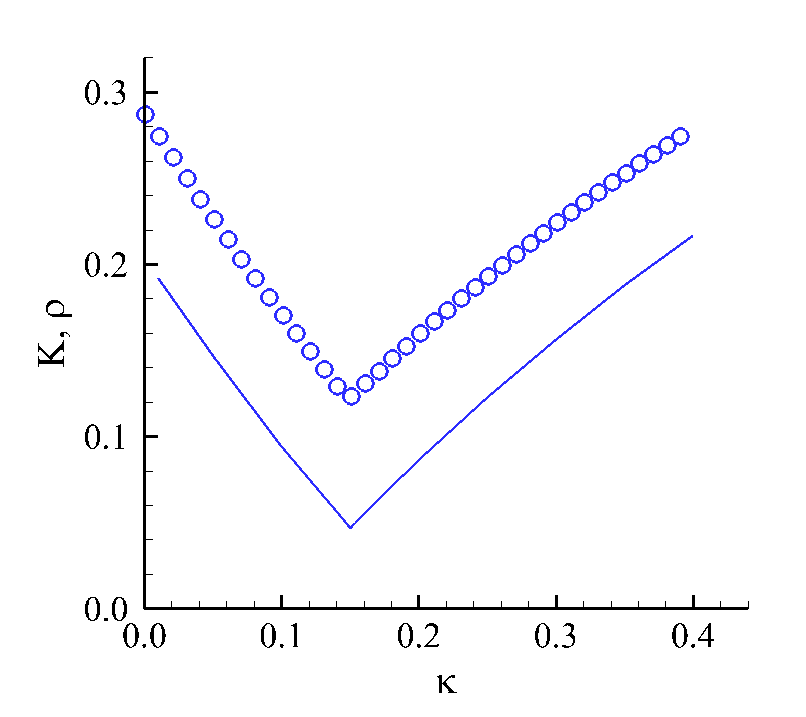}} 
\caption{\small Correlation of $K$ and $\rho $ in the artificial compressible method. a) Case 1. b) Case 2. c) Case 3.} 
\label{artificial_linear_kappa}
\vspace {-0.2 cm}
\end{figure}

Now, let us consider the convergence of artificial compressibility method (\ref {Schwarz_iteration_3}) in association with different values for $\kappa $. 
%Similarly, it is expected that a lower value of  $K$, or, LHS in (\ref {convergence_cond_2}), may correspond to a lower value of $\rho $, thus a faster convergence. Then, 
For above three cases, $K$, together with $\rho $, at different values of $\kappa $  are plotted in Fig. \ref {artificial_linear_kappa}. It is seen that in all of the three cases, $K$ takes the smallest value at $\omega _k=1$. Additionally, $K$ is a good indicator for the value of $\rho $, and both of them take the lowest values at about a same $\kappa $. Therefore, in practical computation, choosing a value of $\kappa $ at the smallest $K$ tends to lead to the smallest $\rho $, which corresponds to the fast convergence speed.   

\begin{figure}[ht] 
\vspace {-0.5 cm}
\centering
\hspace {-0.8 cm}
\subfloat[][]{\label{Jacobi_sp_radius}\includegraphics[width=0.4\linewidth]{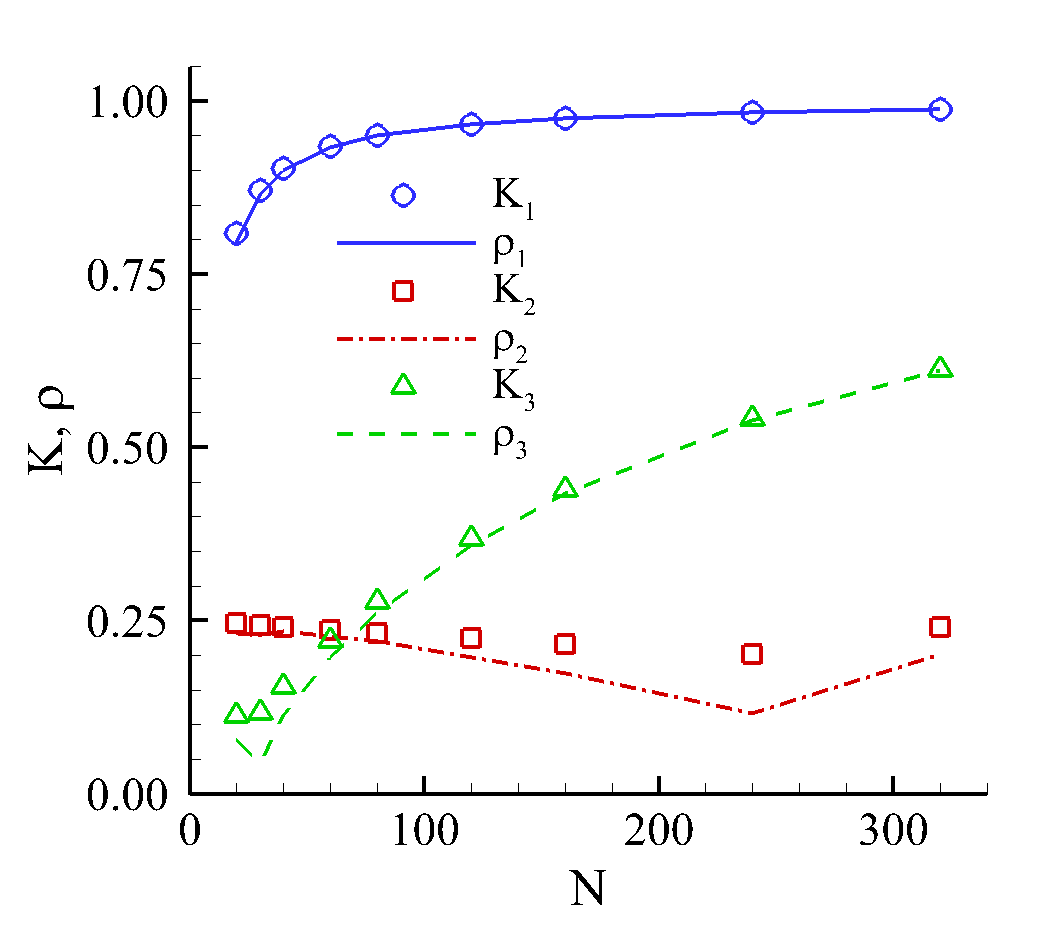}} %\hspace {-0.8 cm}
\subfloat[][]{\label{Pseudo_sp_radius}\includegraphics[width=0.4\linewidth]{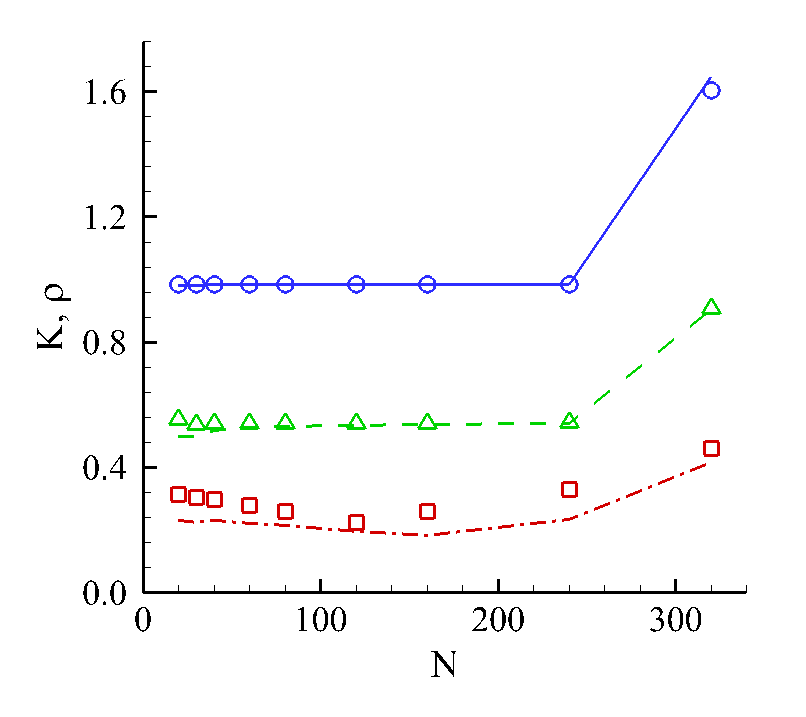}} 
\caption{\small Spectral radius at different mesh resolution. Subscripts for $K$ and $\rho $ indicate the case numbers. $N=I+J-4$, $I$ and $J$ are the number of grid nodes in subdomain 1 and 2, respectively.  $(I,J)=(10,10),(10,20), (20,20),(20,40),(40,40),(40,80),(80,80),(80,160),(160,160)$.
%, corresponding to $\Delta x=0.1176, 0.0741, 0.0541, 0.0351, 0.0260, 0.0171, 0.0127, 0.0084, 0.0063$. 
$\Delta t/\Delta x=0.5$.
%$case1,\eta _1=\eta _2=0.0025; case2, \eta _1=0.1250, \eta _2=0.0625; case3,\eta _1=0.0625, \eta _2=-0.0625$.
a) The Jacobi method. b) The compressibility method. $\kappa =59.25, 0.19$, and $1.19$ for Case 1, 2, and 3, respectively.}
\label{Jacobi_Pseudo_Radius}
\vspace {-0.2 cm}
\end{figure}

\underline {Dependence of convergence speed on grid spacing and time step} Let us look into the convergence speed and start with two identical equations in subdomains. It is seen in above that the iteration matrices for the Jacobi method and the artificial compressibility method are tridiagonal. For an $N\times N$ tridiagonal matrix with elements $p_{-1}$, $p_0$, $p_1$, its eigenvalues read as \cite {Gover:1994}
\begin {equation}
\label{A_radius}
\lambda _i= p_0-2cos\left ( \frac {i\pi }{N+1}\right) \sqrt {p_{-1}p_1}, \ i=1,2, ..., N
\end {equation}
\noindent by which, the spectral radius of the Jacobi method becomes
\begin {equation*}
\rho = 2cos\left ( \dfrac {N\pi }{N+1}\right)\dfrac {\sqrt{\epsilon _k^2-\eta _k^2}}{1+2\epsilon_k+\gamma _k} 
\end {equation*}
\noindent Let $\epsilon _k=const$ and $\epsilon _k\ge \eta _k$, the latter of which is always true when grid spacing is sufficiently small. Then, as $\Delta t, \Delta x$ decrease, the term of ``cos'' increases, and $\eta _k$ and $\gamma _k$ decreases, leading to increase in the term of fraction. As a result, $\rho $ increases monotonically as $\Delta t, \Delta x$ decrease, with an upper limit of $2\epsilon _k/(1+2\epsilon)<1$, and this is 
%what concluded in the literature \cite {Gastaldi1989}, 
seen in Case 1 and 3 in Fig. \ref {Jacobi_sp_radius}. If $\eta _k=const$ and $\epsilon _k< \eta _k$, as the mesh gets finer, the term of fraction does not increase or even decreases, and this may lead to no increase on even a decrease in the spectral radius. A numerical example is Case 2 in the figure. All of these indicate that, in some situations, the spectral radius, or, the convergence speed, is scalable with mesh resolution. Note that the radius may increase again at much finer grids. Similar discussion may be made for the artificial compressibility method, for which the spectral radius becomes 
\begin {equation*} %% kappa = constant? how will the rho change? 0718
\rho =\dfrac{\kappa -(2\epsilon_k+\gamma _k)}{1+\kappa }+\dfrac {2}{1+\kappa }cos\left ( \dfrac {N\pi }{N+1}\right)\sqrt {\epsilon _k^2-\eta _k^2} \\
\end {equation*}
\noindent Again, the spectral radius may be scalable in terms of grid spacing when $\kappa $ takes certain values, and numerical examples are shown in Fig. \ref {Pseudo_sp_radius}. Note that, in studies of waveform relaxation methods, it is concluded that convergence rates deteriorate as grid spacing decrease, e.g., \cite {Gander1998}.  

\underline {Outer and inner iteration} If an iterative scheme is adopted in the subdomains, two types of iteration are involved when marching from time level $n$ to $n+1$: one is the iteration within subdomains, referred to as the inner iteration, and the other is the Schwarz iteration between them, called the outer iteration. Such a situation is commonly encountered in computation of realistic problems, e.g., \cite {Tang2003}. To illustrate the computation, Jacobi iteration (\ref {Schwarz_iteration_2}) is modified as 
\begin{eqnarray}
%\begin {equation}
\label{Schwarz_iteration_4}
\left\{\begin{array}{lll}
(1+2\epsilon  _1+\gamma _1){u_1}_i^l=(\eta _1+\epsilon _1){u_1}_{i-1}^{l-1}+(\epsilon  _1-\eta  _1){u_1}_{i+1}^{l-1}+{u_1}_i^n, & l=1,2, ... L; \ i\le I-1;  \\
{u_1}_I^l={u_2}_2^{m-1} \\
(1+2\epsilon  _2+\gamma _2){u_2}_i^l=(\eta  _2+\epsilon _2){u_2}_{i-1}^{l-1}+(\epsilon  _2-\eta  _2){u_2}_{i+1}^{l-1}+{u_2}_i^n, & l=1,2, ... L; \ i\ge 2; \\
{u_2}_1^l={u_1}_{I-1}^{m-1}
\end{array}\right.  
%\end {equation} 
\end{eqnarray}
\noindent in which $l$ is the inner iteration index, and $m$ is the outer iteration index. After inner iteration  fully converges, ${u_1}_i^l \rightarrow {u_1}_i^m$, ${u_2}_i^l\rightarrow {u_2}_i^m$. Similar formulation can also be made for the artificial compressibility method. A natural question in computation of (\ref {Schwarz_iteration_4}) will be how to arrange the inner and outer iteration to achieve fast convergence, or, after how many times of inner iteration an outer iteration should be started to achieve a minimal count of computational load. 

It is not straightforward to theoretically answer the above question, and a numerical experiment is made to achieve a preliminary understanding. Particularly, computation is made for the three cases in Table \ref {linear_cases} in association with the following initial and boundary conditions
\begin {equation}
g(x)=-sin(\pi x), \  t=0;   \ u_1=0,\ x=-1, \ u_2=0, \ x=1
\label{ic_bc}
\end {equation}
\noindent In the experiment, the corresponding conditions for convergence, i.e., (\ref {convergence_cond_1}) and (\ref {convergence_cond_2}), are satisfied. Also, the computation starts the outer iteration after every 1, 5, and 15 times of the inner iteration, and also after complete convergence of the inner iteration. The results are presented in Table \ref {inner_outer_iteration}. Note that in the computation, for each case, the difference in computational load mainly comes from the numbers of inner iteration, while its number of outer iteration remains about the same. The table shows that the earlier the outer iteration starts, the smaller the amount of computational load. More particularly, the 1:1 strategy is the best; the computation that starts an outer iteration after every inner iteration spends the least amount of computational load.
\begin{table}
\begin{center}
\caption {\small Computational load with different combinations of the inner and outer iteration. `1:1', means after one outer iteration, one inner iteration is made, and `conver.' indicates that after one outer iteration, inner iteration keeps going until convergence. The computational load is defined as the times of computation in each subdomain, or, the total number of times of inner iteration in each subdomain, together with convergence tolerance of $10^{-12}$.}
\begin{tabular}{|c|c|c|c|c|c|c|c|}\hline
Case         & 1:1 & 1:5 & 1:15 & conver.  \\ \hline
1            & 8   & 43  & 102  & 131    \\ \hline
2            & 8   & 43  & 111  & 147    \\ \hline
3            & 8   & 43  & 115  & 152    \\ \hline
\end{tabular}
\vspace {-0.5 cm}
\label{inner_outer_iteration}
\end{center}
\end{table}

\vspace{0.3cm} 

\noindent {\bf 4. Implicit scheme}

\vspace{0.2cm} 

\noindent {\bf 4.1 Convergence speed } 

\vspace{0.2cm} 

An implicit scheme to compute (\ref {adr_dis_1}) in the subdomains reads as 
\begin{equation}
\label{adr_dis_4}
\frac{{u_k}_i^m-{u_k}_i^n}{\Delta t} + a_k\frac{{u_k}_{i+1}^m-{u_k}_{i-1}^m}{2\Delta x}=b_k\frac{{u_k}_{i+1}^m-2{u_k}_i^m+{u_k}_{i-1}^m}{\Delta x^2}-c_k{u_k}_i^m
\end{equation}
\noindent Together with the interface condition, the computation can be formulated as 
\begin{eqnarray}
%\begin {equation}
\label{Schwarz_iteration_5a}
\left\{\begin{array}{lll}
(\eta _1+\epsilon _1){u_1}_{i-1}^m+(1+2\epsilon  _1+\gamma _1){u_1}_i^m+(\epsilon  _1-\eta  _1){u_1}_{i+1}^m+{u_1}_i^n, & i\le I-1;  \
&{u_1}_I^m={u_2}_2^{m-1} \\
(\eta  _2+\epsilon _2){u_2}_{i-1}^m+(1+2\epsilon  _2+\gamma _2){u_2}_i^m+(\epsilon  _2-\eta  _2){u_2}_{i+1}^m+{u_2}_i^n, & i\ge 2; \ &{u_2}_1^m={u_1}_{I-1}^{m-1} 
\end{array}\right.  
%\end {equation} 
\end{eqnarray} 
\noindent The computation in each subdomain needs to solve a linear system. At convergence, ${u_1}^m_i$ and ${u_2}^m_i$ become ${u_1}_i^{n+1}$ and ${u_2}_i^{n+1}$, respectively. In terms of residual, i.e.,  ${e_1}^m_i={u_1}^m_i-{u_1}^{n+1}_i$ and ${e_2}^m_i={u_2}^m_i-{u_2}^{n+1}_i$, above computation becomes
\begin{eqnarray}
\label{Schwarz_iteration_5b}
\left\{\begin{array}{ll}
-(\eta _1+\epsilon _1){e_1}_{i-1}^m+(1+2\epsilon  _1+\gamma _1){e_1}_i^m-(\eta  _1-\epsilon  _1){e_1}_{i+1}^m=0, \ i\le I-1; \ \ &{e_1}^m_I={e_2}^{m-1}_1\\
-(\eta  _2+\epsilon _2){e_2}_{i-1}^m+(1+2\epsilon  _2+\gamma _2){e_2}_i^m-(\eta  _2-\epsilon  _2){e_2}_{i+1}^m=0, \ i\ge 2; \ \ &{e_2}^m_1={e_1}^{m-1}_{I-1} 
\end{array}\right.  
\end{eqnarray}

{\textsc{proposition 4.1}} The contraction factor of computation (\ref {Schwarz_iteration_5b}) is a recursive expression:  
\begin{eqnarray}
\label {convergence_rate_1a}
\bar {\rho }=\dfrac {(\eta _1-\epsilon _1)(\eta _2+\epsilon _2)}{R_1(I-2)R_2(2)}
\end{eqnarray} 
\noindent where $\bar {\rho }={e_1}_i^{m+2}/{e_1}_i^m$, ${e_2}_i^{m+2}/{e_2}_i^m$ and
\begin{eqnarray}
\label {convergence_rate_1b}
\begin{array}{ll}
R_1(i)=1+2\epsilon _1+\gamma _1 +\dfrac {\eta _1^2-\epsilon _1^2}{R_1(i-1)}; \ i=3,...,I-1; \ R_1(2)=1+2\epsilon _1+\gamma _1 \\
R_2(j)=1+2\epsilon _2+\gamma _2 +\dfrac {\eta _2^2-\epsilon _2^2}{R_2(j+1)}; \ j=J-2,...,2; \ R_2(J-1)=1+2\epsilon _2+\gamma _2
\end{array}
\end{eqnarray} 
\noindent in which $I, J\ge 3 $.

{\it Proof}: The proposition can be proved by the induction method. The two equations in (\ref {Schwarz_iteration_5b}) when $i=2 $, $3$ are
\begin{eqnarray*}
\begin{array}{ll}
(1+2\epsilon  _1+\gamma _1){e_1}_2^m+(\eta  _1-\epsilon  _1){e_1}_3^m=0  & \\
-(\eta  _1+\epsilon _1){e_1}_2^m+(1+2\epsilon  _1+\gamma _1){e_1}_3^m+(\eta  _1-\epsilon  _1){e_1}_4^m=0 
\end{array}
\end{eqnarray*}

\noindent In above, without loss of generality (e.g., , $g(x)$ in (\ref {ibvp_adr}) has a finite support), ${e_1}_1^m, {e_2}_3^m=0$ are used. It is readily derived the following by elimination of ${e_1}_2^m$ in above two equations:  
\begin{eqnarray*}
\begin{array}{ll}
(1+2\epsilon  _1+\gamma _1+\dfrac {\eta _1^2-\epsilon _1^2}{R_1(2)}){e_1}_3^m+(\eta  _1-\epsilon  _1){e_1}_4^m=0\\
\end{array}
\end{eqnarray*}

\noindent Then, continue above derivation, for instance, from the above equation and the equation in (\ref {Schwarz_iteration_5b}) when $i=4$ in subdomain 1. It is expected that  
\begin{eqnarray*}
\begin{array}{ll}
\left (1+2\epsilon _1+\gamma _1+\dfrac{\eta _1^2-\epsilon _1^2}{R_1(I-3)}\right ){e_1}_{I-2}^m+(\eta  _1-\epsilon  _1){e_1}_{I-1}^m=0\\  
\end{array}
\end{eqnarray*}

\noindent The equation (\ref {Schwarz_iteration_5b}) in subdomain 1 when $i=I-1$ is 
\begin{eqnarray*}
\begin{array}{ll}
-(\eta  _1+\epsilon _1){e_1}_{I-2}^m+(1+2\epsilon  _1+\gamma _1){e_1}_{I-1}^m+(\eta  _1-\epsilon  _1){e_1}_I^m=0 \\ 
\end{array}
\end{eqnarray*}

\noindent Using the two equations and noticing that the first parenthesis in the former is actually $R_1(I-2)$, it is readily derived that 
\begin{eqnarray*}
\begin{array}{ll}
\left (1+2\epsilon _1+\gamma _1+\dfrac{\eta _1^2-\epsilon _1^2}{R_1(I-2)}\right ){e_1}_{I-1}^m+(\eta _1-\epsilon  _1){e_1}_I^m=0\\  
\end{array}
\end{eqnarray*} 

\noindent which, after substitution of the interface conditions in (\ref {Schwarz_iteration_5b}) and change of index $m$ in the latter, becomes 
\begin{eqnarray*}
\begin{array}{ll}
\left (1+2\epsilon _1+\gamma _1+\dfrac{\eta _1^2-\epsilon _1^2}{R_1(I-2)}\right ){e_2}_1^{m+2}+(\eta _1-\epsilon  _1){e_2}_2^m=0\\  
\end{array}
\end{eqnarray*} 

\noindent Similarly, starting from the two equations in subdomain 2 when $i=J-1, J-2$, one can derived that 
\begin{eqnarray*}
\begin{array}{ll}
-(\eta  _2+\epsilon _2){e_2}_1^m+\left (1+2\epsilon  _2+\gamma _2+\dfrac{\eta _2^2-\epsilon _2^2}{R_2(2)}\right ){e_2}_2^m=0  
\end{array}
\end{eqnarray*}

\noindent Above two equations results in (\ref {convergence_rate_1a}) when $\bar {\rho }={e_2}_1^{m+2}/{e_2}_1^m$. Moreover, from the latter, it is readily check that ${e_2}_2^{m+2}/{e_2}_2^m={e_2}_1^{m+2}/{e_2}_1^m$. Also, as seen in above proof, the relation like the latter equation holds at $i>2$, too. These indicate that (\ref {convergence_rate_1a}) is true at all $i$ in subdomain 1. In the same way, one can prove above for subdomain 1. These completes the proof. 

{\it Remark 4.1} Since the contraction factor is the same at all nodes, it is readily seen that such factor is actually the global convergence speed, e.g., $\bar {\rho }=||{e_1}_i^{m+2}||_2/||{e_1}_i^m||_2$.
% ||{e_2}_i^m||_2/||{e_2}_i^{m-2}||_2$.  

In order to validate the derived convergence rate (\ref {convergence_rate_1a}), a numerical experiment is made in computation of the three cases in Table \ref {linear_cases} in conjunction with initial and boundary conditions (\ref {ic_bc}). The theoretical contraction factor from (\ref {convergence_rate_1a}) and numerical one computed directly from the solutions of the experiment are listed in Table \ref {convergence_rate}, which shows that the two factors are basically identical, validating the derived one.  

\begin{table}[ht]
  \centering
  \caption{\small Theoretical and numerical convergence speeds. $(I,J)=(40,80)$, $\Delta t, \Delta x=0.0171$.} 
  \label{convergence_rate}
  \begin{tabular}{|c|c|c|c|c|c|}
\hline
      Cases   &  $\bar {\rho }_{theo}$ &  $\bar {\rho }_{num}$   \\ \hline
        1     &   0.691235   &   0.691197   \\ \hline
        2     &  -0.017007   &  -0.017007   \\ \hline
        3     &   0.082962   &   0.082962   \\ \hline
  \end{tabular}
\end{table}

\vspace{0.8cm} 

\noindent {\bf 4.2 Optimal interface condition}

\vspace{0.2cm} 

In order to speed up convergence in computation of (\ref {Schwarz_iteration_5a}), the following optimal interface condition is adopted: 
\begin {equation}
\label{optimal_int}
\begin{array}{ll}
({u_1}^m_I-{u_1}^m_{I-1})+\alpha {u_1}^m_I=({u_2}^{m-1}_2-{u_2}^{m-1}_1)+\alpha {u_2}^{m-1}_2 \\
({u_2}^m_2-{u_2}^m_1)+\beta {u_2}^m_1=({u_1}^{m-1}_I-{u_1}^{m-1}_{I-1})+\beta {u_1}^{m-1}_{I-1} 
\end{array}
\end {equation} 
\noindent in which $\alpha , \beta =consts$. Such optimal interface condition was first proposed for waveform relaxation of circuit problems by Gander et al., and it was effective in speeding up the convergence  \cite {Gander:2002}. They show that, if $(\alpha +1)(\beta -1)+1\ne 0$, the optimal interface condition recovers to the classic condition, i.e., the interface condition in (\ref {Schwarz_iteration_5a}), at convergence as $m \rightarrow \infty$. Also, obviously, the optimal condition returns to the classic condition as $\alpha , \beta \rightarrow \infty$. 

With the optimal interface condition, (\ref {Schwarz_iteration_5b}) is modified as
\begin{eqnarray}
%\begin {equation}
\label{Schwarz_iteration_6}
\left\{\begin{array}{ll}
-(\eta _1+\epsilon _1){e_1}_{i-1}^m+(1+2\epsilon  _1+\gamma _1){e_1}_i^m-(\eta _1-\epsilon  _1){e_1}_{i+1}^m=0, \ i\le I-1 \\
\ \ ({e_1}^m_I-{e_1}^m_{I-1})+\alpha {e_1}^m_I=({e_2}^{m-1}_2-{e_2}^{m-1}_1)+\alpha {e_2}^{m-1}_1\\

-(\eta  _2+\epsilon _2){e_2}_{i-1}^m+(1+2\epsilon  _2+\gamma _2){e_2}_i^m-(\eta  _2-\epsilon  _2){e_2}_{i+1}^m=0, \ i\ge 2 \\ 
\ \ ({e_2}^m_2-{e_2}^m_1)+\beta {e_2}^m_1=({e_1}^{m-1}_I-{e_1}^{m-1}_{I-1})+\beta {e_1}^{m-1}_I 
\end{array}\right.  
%\end {equation} 
\end{eqnarray}
\noindent The convergence speed for computation of (\ref {Schwarz_iteration_6}) can be derived by following the same method and steps adopted to those in computation of (\ref {Schwarz_iteration_5b}), and this is based on a fact that the grid nodes involved in the former are exactly same as those in the latter. The derived results are summarized as follows.   

{\textsc{proposition 4.2}} The contraction factor in computation of (\ref {Schwarz_iteration_6}) is   
\begin{eqnarray}
\label{convergence_rate_2}
\bar {\rho }=-\frac {(\alpha +1)(\eta _2 +\epsilon _2)-R_2(2)}{(\alpha +1)R_1(I-1)+\eta _1-\epsilon _1}
            \cdot \frac {(\beta -1)(\eta _1-\epsilon _1)-R_1(I-1)}{(\beta -1)R_2(2)+\eta _2 +\epsilon _2}
\end{eqnarray} 
\noindent in which $I,J \ge 3$. 

In expression (\ref {convergence_rate_2}), letting the numerators be zero while keeping denominators be non-zero leads to zero contraction factor, i.e., $\bar {\rho }=0$, and the following optimal $\alpha $ and $\beta $: 
%\begin{eqnarray}
\begin {equation}
\label{optimal_1}
\alpha =\dfrac {R_2(2)}{\eta _2 +\epsilon _2}-1, \ \
\beta  =\dfrac {R_1(I-1)}{\eta _1 -\epsilon _1}+1 
\end {equation} 
%\end{eqnarray} 
\noindent together with condition $R_1(I-1)R_2(2)+(\eta _1 -\epsilon _1 )(\eta _2 +\epsilon _2)\ne 0$. From (\ref {convergence_rate_2}) and (\ref {optimal_1}), it is derived that
\begin{eqnarray*}
\label{opt_linear_1st_derivative}
\begin{array}{lll}
\dfrac {\partial \bar {\rho }}{\partial \alpha }
=-\dfrac {R_1(I-1)R_2(2)+(\eta _1 -\epsilon _1 )(\eta _2 +\epsilon _2)}{[(\alpha +1)R_1(I-1)+(\eta _1-\epsilon _1)]^2}
 \cdot \dfrac {(\beta -1)(\eta _1-\epsilon _1)-R_1(I-1)}{(\beta -1)R_2(2)+(\eta _2 +\epsilon _2)} \\
\dfrac {\partial \bar {\rho }}{\partial \beta }
=-\dfrac {(\alpha +1)(\eta _2 +\epsilon _2)-R_2(2)}{(\alpha +1)R_1(I-1)+(\eta _1-\epsilon _1)}
            \cdot \dfrac {R_1(I-1)R_2(2)+(\eta _1 -\epsilon _1 )(\eta _2 +\epsilon _2)}{[(\beta -1)R_2(2)+(\eta _2 +\epsilon _2)]^2}
\end{array}
\end{eqnarray*}
\noindent additionally, 
\begin{eqnarray*}
\label{opt_linear_2nd_derivative}
\begin{array}{lll}
\dfrac {\partial ^2\bar {\rho }}{\partial \alpha ^2}
=-\dfrac {2R_1(I-1)[R_1(I-1)R_2(2)+(\eta _1 -\epsilon _1 )(\eta _2 +\epsilon _2)]}{[(\alpha +1)R_1(I-1)+(\eta _1-\epsilon _1)]^3}
 \cdot \dfrac {(\beta -1)(\eta _1-\epsilon _1)-R_1(I-1)}{(\beta -1)R_2(2)+(\eta _2 +\epsilon _2)} \\
\dfrac {\partial ^2\bar {\rho }}{\partial \beta ^2}
=-\dfrac {(\alpha +1)(\eta _2 +\epsilon _2)-R_2(2)}{(\alpha +1)R_1(I-1)+(\eta _1-\epsilon _1)}
 \cdot \dfrac {2R_2(2)[R_1(I-1)R_2(2)+(\eta _1 -\epsilon _1 )(\eta _2 +\epsilon _2)]}{[(\beta -1)R_2(2)+(\eta _2 +\epsilon _2)]^3} \\
\dfrac {\partial ^2\bar {\rho }}{\partial \alpha \partial \beta }
=-\dfrac {[R_1(I-1)R_2(2)+(\eta _1 -\epsilon _1 )(\eta _2 +\epsilon _2)]^2}{[(\alpha +1)R_1(I-1)+(\eta _1-\epsilon _1)]^2 [(\beta -1)R_2(2)+(\eta _2 +\epsilon _2)]^2}
\end{array}
\end{eqnarray*}
\noindent Under conditions (\ref {optimal_1}), 
%and $ \eta _1\ne \epsilon _1$, which is always true as $\Delta t$ and $\Delta x$ are sufficiently small, 
it is readily verified that 
\begin {equation*}
\label{optimal_convergence_rate}
\begin{array}{lll}
\dfrac {\partial \bar {\rho }}{\partial \alpha}= 0, \ \ 
\dfrac {\partial \bar {\rho }}{\partial \beta}= 0 \\
\dfrac {\partial ^2\bar {\rho }}{\partial \alpha ^2}= 0, \ \ 
\dfrac {\partial ^2\bar {\rho }}{\partial \beta ^2}= 0,  \ \ 
\dfrac {\partial ^2\bar {\rho }}{\partial \alpha \partial \beta }\ne 0  
\end{array}
\end {equation*} 
\noindent As a result, one has 
\begin {equation*}
\frac {\partial ^2\bar {\rho }}{\partial \alpha ^2}
\dfrac {\partial ^2\bar {\rho }}{\partial \beta ^2}
-\left (\frac {\partial ^2\bar {\rho }}{\partial \alpha \partial \beta }\right )^2<0
\end {equation*} 
%{\bf it should be 2nd derivatives in the above first term ?}
\noindent With above, it is known that ($\alpha $, $\beta $) given by (\ref {optimal_1}) is a saddle point of $\bar \rho $ \cite {Stewart2005}. All above proves the following theorem. 

{\textsc{Theorem 4.1}} At the optimal value of ($\alpha $, $\beta $) given in (\ref {optimal_1}), the contraction factor (\ref {convergence_rate_2}) is zero. Moreover,  the optimal value of ($\alpha $, $\beta $) is a saddle point of $\bar \rho $. 

{\it Remark 4.2} Actually, contraction factor (\ref {convergence_rate_2}) becomes zero when either $\alpha $ takes the value in (\ref {optimal_1}), regardless that of $\beta $, or, $\beta $ takes the value in (\ref {optimal_1}), regardless that of $\alpha $. Also, the zero contraction factor occurs only when ($\alpha ,\beta $) are given as (\ref {optimal_1}). 

{\it Remark 4.3} The past work is at the continuous or semi-discretization level, and its optimal values for ($\alpha ,\beta $) can only be derived in a transformed space. As a result, search and approximation for them at the algebraic level have to be made in actual computation, e.g., \cite {Gander2004}. As a distinct feature in this study, the optimal ($\alpha , \beta $) can be directly computed by (\ref {optimal_1}), without any approximation. Therefore, the derived theoretical convergence speed exactly matches that in practical computation, as shown in the following numerical examples.  

\begin{figure}[ht] 
\vspace {-0.5 cm}
\centering
\hspace {-0.8 cm}
\subfloat[][]{\label{linear_opt_convergence}\includegraphics[width=0.4\linewidth]{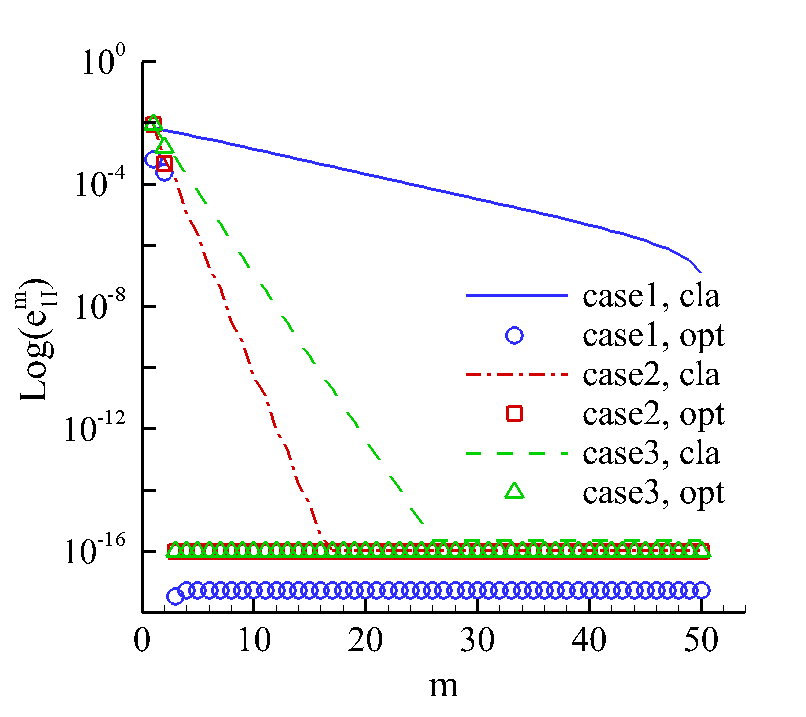}} %\hspace {-0.8 cm}
\subfloat[][]{\label{saddle_pt}\includegraphics[width=0.425\linewidth]{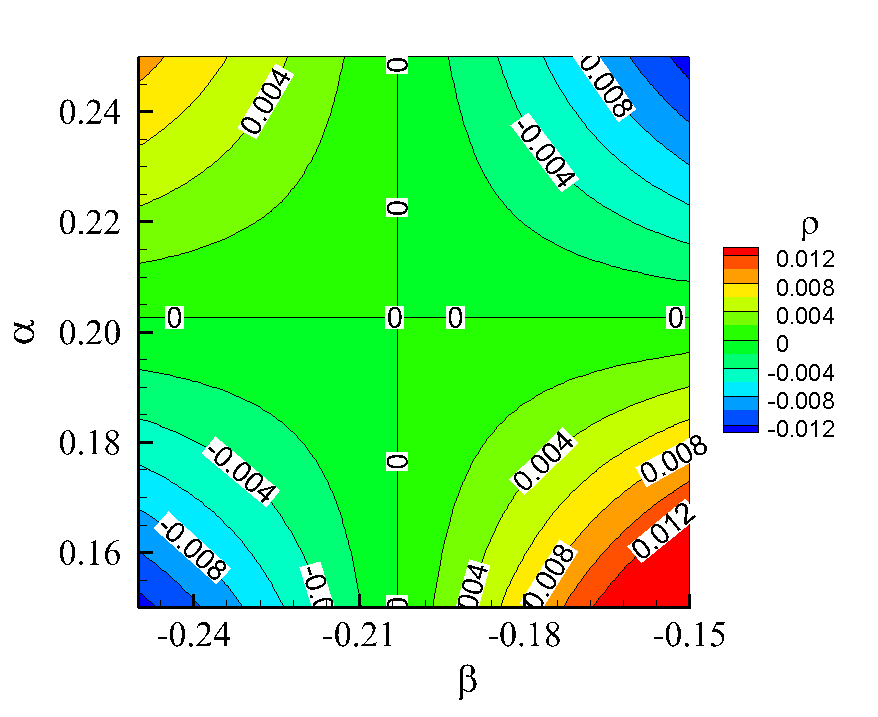}} 
\caption{\small  Computation with optimal interface condition. a) Convergence processes. b) Contraction factor. Case 1, (I,J)=(40,80).
%$\alpha $ and $\beta $, with optimal values ($\alpha, \beta $)=($2.02577\times 10^{-1},-2.02988\times 10^{-1}$).
} 
\label{Convergence_and_saddle_pt}
\vspace {-0.2 cm}
\end{figure}

Again, a numerical experiment is made on the three cases in Table \ref {linear_cases} with initial and boundary conditions (\ref {ic_bc}), and the results are shown in Table \ref {linear_opt_convergence_rate}. The table shows that, at the optimal values for ($\alpha ,\beta $), the numerical contraction factors in the numerical solutions are very small or almost zeros. These factors are not exactly zero, which is the theoretical value, and this is attributed to the lack of enough accuracy in the computer. Calculations indicate that the value of $\bar {\rho }$ is subtle to the values of $\alpha $ and $\beta $, or, enough digits should be kept in calculation of (\ref {optimal_1}) to achieve zero for a contraction factor. To directly illustrate the effectiveness of the optimal interface condition, the convergence residuals in the experiment are shown in Fig. \ref {linear_opt_convergence}. The figure shows that, in comparison to the classic condition, the optimal condition greatly speeds up the convergence. In view that the contraction factor is the ratio of residual between two subsequent times of iteration, a zero convergence factor implies that convergence will be achieved with no more than two times of convergence. The figure shows that indeed the convergence is achieved in about two times of iteration. Moreover, Fig. \ref {saddle_pt} illustrates a distribution of the contraction factor on the $\alpha $-$\beta $ plane, and it clearly shows that the values of $\alpha $ and $\beta $ given by (\ref {optimal_1}) is a saddle point. 

\begin{table}[ht]
  \centering
  \caption{\small Optimal ($\alpha ,\beta $) and convergence speed. (I,J)=(40,80), $\Delta t, \Delta x=0.017$.} 
  \label{linear_opt_convergence_rate}
  \begin{tabular}{|c|c|c|c|}
\hline
Cases & $\alpha$, $\qquad \qquad $ $\beta $                & $\bar {\rho }_{num}$   \\ \hline
1     & 2.02577$\times 10^{-1}$, -2.02988$\times 10^{-1}$  & 1.31$\times 10^{-14}$  \\ \hline
2     & 5.14615$\times 10^{0}$,  1.05666$\times 10^{1}$    & -6.12$\times 10^{-15}$  \\ \hline
3     & 1.74955$\times 10^{0}$, -3.38389$\times 10^{0}$    & 0                  \\ \hline
  \end{tabular}
\end{table}

\vspace{0.3cm} 

\noindent {\bf 5. Extension to Burgers' equation}

\vspace{0.2cm} 

\noindent {\bf 5.1 Explicit scheme}

\vspace{0.2cm} 

Consider the coupling between two viscous Burgers' equations with reaction terms. Now the discretization (\ref {adr_dis_1}) becomes
\begin{equation}
\label{adr_dis_5}
\frac{{u_k}_i^{n+1}-{u_k}_i^n}{\Delta t} + {u_k}_i^{n+1}\frac{{u_k}_{i+1}^{n+1}-{u_k}_{i-1}^{n+1}}{2\Delta x}=b_k\frac{{u_k}_{i+1}^{n+1}-2{u_k}_i^{n+1}+{u_k}_{i-1}^{n+1}}{\Delta x^2}-c_k{u_k}_i^{n+1},
\end{equation} 
\noindent which can be computed with an extension of (\ref {adr_dis_2}), i.e., the Jacobi method for the linear advection-diffusion-reaction equation, particularly, by letting its $a_k$ be replaced by ${u_k}_i^{m-1}$. With such an extension, computation of above discretization in association with the classic interface condition leads to an explicit algorithm 
\begin{eqnarray}
\label{Schwarz_iteration_7}
\left\{\begin{array}{lll}
(1+2\epsilon  _1+\gamma _1){u_1}_i^m=({\eta _1}_i^{m-1}+\epsilon _1){u_1}_{i-1}^{m-1}+(\epsilon  _1-{\eta  _1}_i^{m-1}){u_1}_{i+1}^{m-1}+{u_1}_i^n, & i\le I-1;  \\
\qquad {u_1}_I^m={u_2}_2^m \\
(1+2\epsilon  _2+\gamma _2){u_2}_i^m=({\eta  _2}_i^{m-1}+\epsilon _2){u_2}_{i-1}^{m-1}+(\epsilon  _2-{\eta  _2}_i^{m-1}){u_2}_{i+1}^{m-1}+{u_2}_i^n, & i\ge 2; \\ 
\qquad {u_2}_1^m={u_1}_{I-1}^m 
\end{array}\right.  
\end{eqnarray}
\noindent A difference of above iteration from the previous iteration is that now the value of $\eta $ becomes solution dependent. In order to analyze above iteration, linearization is made by replacing its ${\eta _k}_i^{m-1}$ with ${\eta _k}_i^n=\Delta t{u_k}_i^n/(2\Delta x)$ ($k=1,2$). As a result, the iteration between the two time levels can be formulated in form of (\ref {Schwarz_iteration_2}) with $\eta _k$ be replaced by ${\eta  _k}_i^n$ , and the associated coefficient matrix $A$ remains constant during the iteration while marching from time level $n$ to $n+1$. According to the Scarborough criterion, the following conclusion is achieved (Note that LHS is not a constant but changes with $i,n$).   

{\textsc{proposition 5.1}} A sufficient condition for the linearized version of (\ref {Schwarz_iteration_7}) to converge is
\begin {equation}
\label{convergence_cond_3}
\max\limits_{i, k}\bigg\{
\frac {|\epsilon_k+{\eta_k }^n_i|+|{\eta_k }_i^n-\epsilon_k|}{|1+2\epsilon_k+\gamma _k|}\bigg \}< 1    
\end {equation} 

Discretization (\ref {adr_dis_5}) can also be computed using the artificial compressibility method, and this is realized by replacing $a_k$ with ${u_k}_i^{m-1}$ in (\ref {adr_dis_3}). Consequently, method (\ref {Schwarz_iteration_3}) is extended as
\begin{eqnarray}
\label{Schwarz_iteration_8}
\left\{\begin{array}{lll}
(1+\kappa ){u_1}_i^m
=({\eta _1}_i^{m-1}+\epsilon _1){u_1}_{i-1}^{m-1}
+\left (\kappa -2\epsilon  _1-\gamma _1\right ){u_1}_i^{m-1}
+(\epsilon  _1-{\eta  _1}_i^{m-1}){u_1}_{i+1}^{m-1}+{u_1}_i^n, \\ 
\qquad i\le I-1; \ {u_1}_I^m={u_2}_2^m \\
(1+\kappa ){u_2}_i^m
=({\eta _2}_i^{m-1}+\epsilon _2){u_2}_{i-1}^{m-1}
+\left (\kappa -2\epsilon  _2-\gamma _2\right ){u_2}_i^{m-1}
+(\epsilon  _2-{\eta  _1}_i^{m-1}){u_2}_{i+1}^{m-2}+{u_2}_i^n, \\ 
\qquad i\ge 2;  \ {u_2}_1^m={u_2}_{I-1}^m
\end{array}\right.  
\end{eqnarray}
\noindent By linearization to replace ${\eta _k}_i^{m-1}$ with ${\eta _k}_i^n$ ($k=1,2$) and with a similar discussion, the following is obtained. 

{\textsc{proposition 5.2}} A sufficient condition for the linearized version of (\ref {Schwarz_iteration_8}) to converge is 
\begin {equation}
\label{convergence_cond_4}
\max\limits_{i, k}\bigg\{
\bigg| \dfrac {\epsilon_k+{\eta_1 }^n_i}{1+\kappa } \bigg|+\bigg| \dfrac{\kappa -(2\epsilon_k+\gamma _k)}{1+\kappa }\bigg|+\bigg|\dfrac {\epsilon_k-{\eta_k }^n_i}{1+\kappa }\bigg|\bigg \}
<1
\end {equation}

The analysis for computation of (\ref {Schwarz_iteration_7}) and (\ref {Schwarz_iteration_8}) can be made similarly as in Section 3.3 for the linear advection-diffusion-reaction equations. As an illustration of the analysis, numerical experiments are made on four cases in Table \ref {Burgers_cases}. Case 4 and 5 are respectively diffusion and advection dominated, respectively, and Case 6 and 7 are a combination of them. Initial and boundary conditions in (\ref {ic_bc}) are used in Case 4 and 7. For case 5 and 6, the initial condition becomes piecewise: $g(x)= 1, x\in [-0.5,0.5], \ g(x)=0, x\in [-1,-0.5), (0.5,1]$. 

\begin{table}[ht]
  \centering
  \caption{\small Cases for numerical experiments on the viscous Burgers equations. } 
  \label{Burgers_cases}
  \begin{tabular}{|c|c|c|c|c|c|c|}
\hline
       Cases & $b_1$ & $c_1$ & $b_2$ & $c_2$ \\ \hline
        4    & 0.5  & 0.0   & 0.5    & 0.0   \\ \hline       
        5    & 0.02 & 2.0   & 0.02   & 3.0   \\ \hline                  
        6    & 0.1  & 0.2   & 0.2    & 0.1   \\ \hline 
        7    & 0.5  & 0.0   & 0.51   & 0.0   \\ \hline 
  \end{tabular}
\end{table}

\begin{figure}[t] 
\centering
\subfloat[][]{\label{t1}\includegraphics[width=0.4\linewidth]{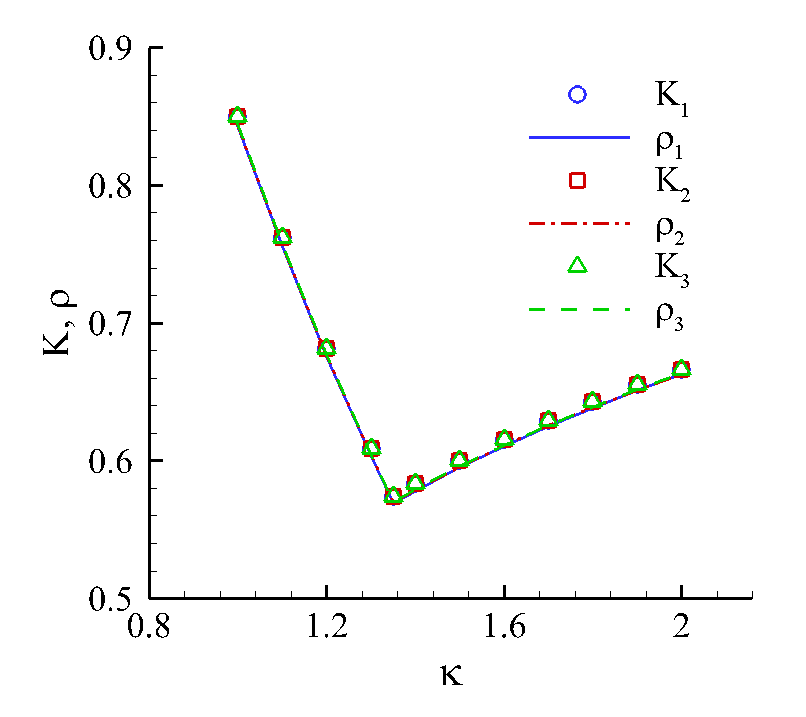}} % \hspace {-0.8 cm}
\subfloat[][]{\label{t2}\includegraphics[width=0.4\linewidth]{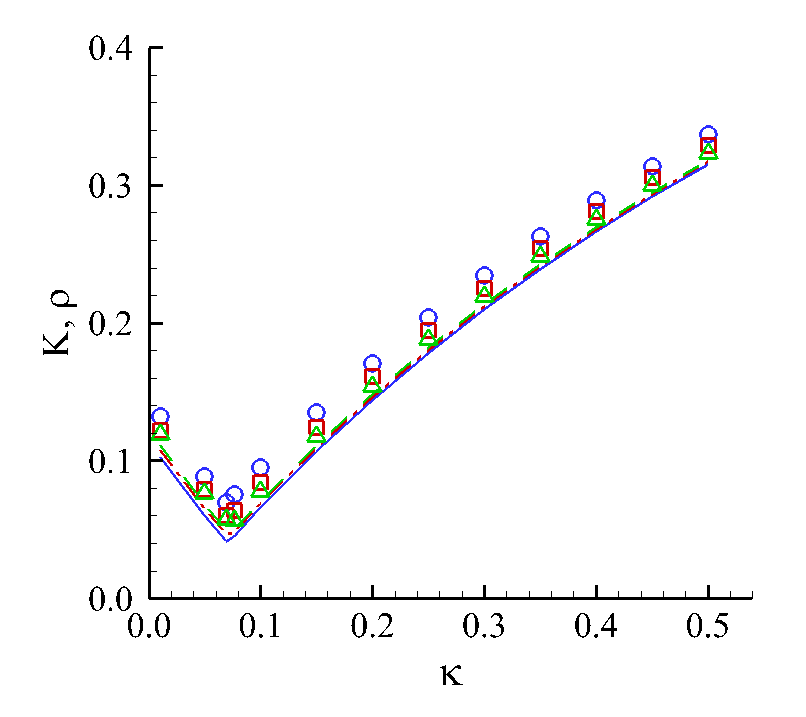}} \\
\subfloat[][]{\label{t3}\includegraphics[width=0.4\linewidth]{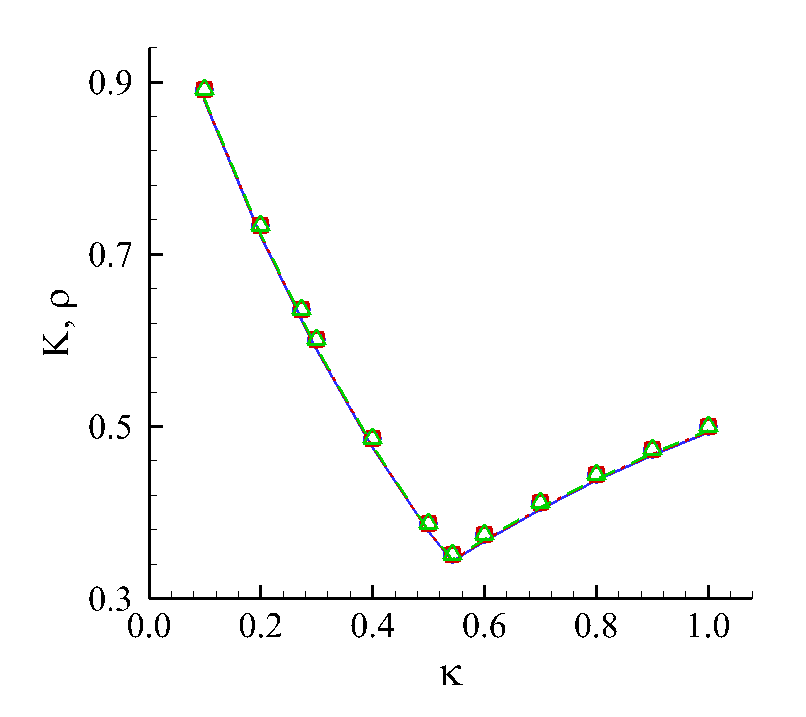}}
\subfloat[][]{\label{t4}\includegraphics[width=0.4\linewidth]{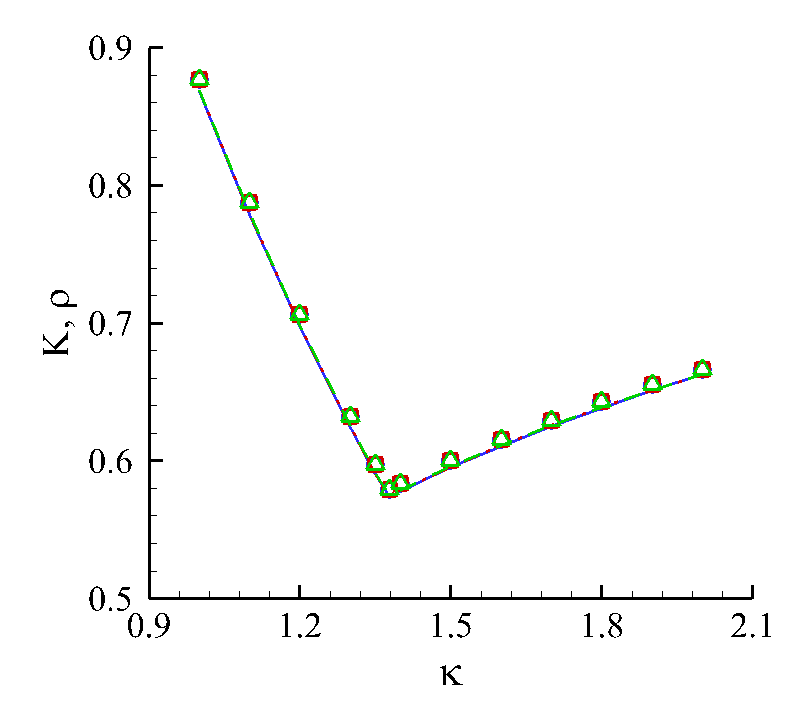}}
\caption{\small Correlation of $K$ and $\rho $ in the artificial compressible method for the Burgers equations. Subscript 1, 2, and 3 indicate t=0.074, 0.222, and 0.370, respectively. a) Case 4. b) Case 5. c) Case 6. d) Case 7.}
\label{artificial_nonlinear_kappa}
\vspace {-0.2 cm}
\end{figure}

Table \ref {SC_radius_Burgers} presents the values for LHS of (\ref {convergence_cond_3}) and (\ref {convergence_cond_4}), together with those for the spectral radii. It is seen that when these two conditions are satisfied, all radii are less than 1, and thus convergence in these linearized versions of (\ref {Schwarz_iteration_7}) and (\ref {Schwarz_iteration_8}) is guaranteed. Noted that since the equations become nonlinear, particularly, in general,  ${\eta _k}_i^n$ changes with solutions, and its values are distinct at different time steps. To illustrate this, results at three moments in time are presented in the table. In computation of (\ref {Schwarz_iteration_8}), effects of $\kappa $ on $K$ and $\rho$ are plotted in Fig. \ref {artificial_nonlinear_kappa}. Again, it is seen that the value of $\kappa $ corresponding to the lowest $K$ is a good approximation of its value for the lowest $\rho $, which results in the fastest convergence. 
%As $\kappa $ takes such a value, the LHS in (\ref {convergence_cond_4}) and spectral radius are also presented in Table \ref {SC_radius_Burgers}. 

\begin{table}[ht]
\centering
\caption{\small LHS of (\ref {convergence_cond_3}) and (\ref {convergence_cond_4}), and corresponding spectral radius. (I,J)=(10,20), $\Delta t, \Delta x=0.0741$. For the artificial compressibility method, $\kappa=1.350$, $0.069$,  $0.542$, and $1.378 $ are used in Case 4, 5, 6, and 7, respectively. These values of $\kappa $ correspond to the lowest values of $K$, as shown in Fig. \ref {artificial_nonlinear_kappa}.
} 
\begin{tabular}{|c|c|c|c|c|c|c|}
\hline
 Time & \multicolumn{2}{c|}{ 0.074 } & \multicolumn{2}{c|}{ 0.222 } &
\multicolumn{2}{c|}{ 0.370 } \\ \hline
 Method, Case &  K  &  $\rho$  & K & $\rho$ & K & $\rho$  \\ \hline
 Jacobi, 4 & 0.574468 & 0.558299 & 0.574468 & 0.566377 & 0.574468 & 0.568848  \\ \hline
 Jacobi, 5 & 0.080057 & 0.055298 & 0.057561 & 0.041587 & 0.050176 & 0.044854  \\ \hline
 Jacobi, 6 & 0.350481 & 0.340231 & 0.350481 & 0.342112 & 0.350481 & 0.343251  \\ \hline
 Jacobi, 7 & 0.574468 & 0.572429 & 0.574468 & 0.573695 & 0.574468 & 0.574239  \\ \hline
Artificial, 4  & 0.574468 & 0.568429 & 0.574468 & 0.569622 & 0.574468 & 0.570135  \\ \hline
Artificial, 5  & 0.081437 & 0.054665 & 0.058767 & 0.047987 & 0.051325 & 0.041229  \\ \hline
Artificial, 6  & 0.635802 & 0.623317 & 0.635802 & 0.625591 & 0.635802 & 0.626952  \\ \hline
Artificial, 7  & 0.597447 & 0.589117 & 0.597447 & 0.590465 & 0.597447 & 0.591051  \\ \hline
\end{tabular}
\label{SC_radius_Burgers}
\end{table}

Lastly, let us look into the issue of inner and outer iteration, and extend Jacobi iteration (\ref {Schwarz_iteration_4}) by replacing $\eta _k$ with ${\eta _k}_i^{l-1}$. The results of computation for the cases in Table \ref {Burgers_cases} are presented in Table \ref {Burgers_inner_outer_iteration}, and they show a same trend observed in the linear advection-diffusion-reaction equations. Particularly, the earlier to start the outer iteration, the smaller the computational load, or, the quicker the overall convergence. The smallest computational load occurs at the 1:1 strategy. Therefore, an understanding is that the 1:1 strategy leads to the fastest convergence in the examples of this paper, regardless the equations are linear or nonlinear. Such a strategy has been used in practical problems \cite {Tang2003}. 

\begin{table}[ht]
\caption {Computational load associated with different combinations of the inner and outer iteration, for the Burgers equations, together with convergence tolerance of $10^{-12}$, at $t=0.170$.}
\begin{center}
\begin{tabular}{|c|c|c|c|c|c|c|c|c|}\hline
Case No. & 1:1 & 1:5 & 1:15 & conver.  \\ \hline
4  & 7  & 63  & 72  & 120  \\ \hline
5  & 9  & 70  & 86  & 132   \\ \hline
6  & 8  & 66  & 78  & 125   \\ \hline
7  & 7  & 62  & 72  & 120  \\ \hline
\end{tabular}
\label{Burgers_inner_outer_iteration}
\end{center}
\end{table}

\vspace{0.3cm} 

\noindent {\bf 5.2 Implicit scheme}

\vspace{0.2cm} 

Now consider computation of discretization (\ref {adr_dis_5}) by an implicit scheme, and this is realized by replacing its $a_k$ with ${u_k}_i^{m-1}$ in scheme (\ref{adr_dis_4}). With such modification, the computational problem becomes
\begin{eqnarray}
\label{Schwarz_iteration_10}
\left\{\begin{array}{lll}
({\eta _1}_i^{m-1}+\epsilon _1){u_1}_{i-1}^m+(1+2\epsilon  _1+\gamma _1){u_1}_i^m+(\epsilon  _1-{\eta  _1}_i^{m-1}){u_1}_{i+1}^m+{u_1}_i^n, \ i\le I-1;  \\
\qquad {u_1}_I^m={u_2}_2^{m-1} \\
({\eta  _2}_i^{m-1}+\epsilon _2){u_2}_{i-1}^m+(1+2\epsilon  _2+\gamma _2){u_2}_i^m+(\epsilon  _2-{\eta  _2}_i^{m-1}){u_2}_{i+1}^m+{u_2}_i^n, \  i\ge 2; \\
\qquad {u_2}_1^m={u_1}_{I-1}^{m-1} 
\end{array}\right.  
\end{eqnarray} 
\noindent which requires to solve linear systems in subdoamins at each iteration over $m$.

Let us analyze computation of (\ref {Schwarz_iteration_10}) by linearization with ${\eta _1}_i^{m-1}$ and $ {\eta _2}_i^{m-1}$ being replaced by ${\eta _1}_i^n$ and ${\eta _2}_i^n$, respectively. As a result, (\ref {Schwarz_iteration_10}) becomes a linear iterative scheme. Follow the proof of Proposition 4.1, still a contraction factor can be derived on the basis of linearization as follows. 

{\textsc{proposition 5.3}} The contraction factor for computing the linearized Burgers equations, (\ref {Schwarz_iteration_10}), is
\begin{eqnarray}
\label{convergence_rate_3a}
\bar {\rho }=\dfrac {({\eta _1}_{I-1}^n-\epsilon _1)({\eta _2}_2^n+\epsilon _2)}{R_1(I-1)R_2(2)}
\end{eqnarray} 
\noindent where $ \bar {\rho }={e_1}_i^{m+2}/{e_1}_i^m$, $ \bar {\rho }={e_2}_i^{m+2}/{e_2}_i^m$, and 
\begin{eqnarray}
\label{convergence_rate_3b}
\begin{array}{ll}
R_1(i)=1+2\epsilon _1+\gamma _1 +\dfrac {({\eta _1}_{i-1}^n-\epsilon _1)({\eta _1}_i^n+\epsilon _1)}{R_1(i-1)};\ i=3,...,I-1,  \ R_1(2)=1+2\epsilon _1+\gamma _1  \\
R_2(j)=1+2\epsilon _2+\gamma _2 +\dfrac {({\eta _2}_j^n-\epsilon _2)({\eta _2}_{j+1}^n+\epsilon _2)}{R_2(j+1)};\ j=J-2,...,2, \ R_2(J-1)=1+2\epsilon _2+\gamma _2
\end{array}
\end{eqnarray} 
\noindent in which $I,J\ge 3$.

From (\ref {convergence_rate_3a}), it is seen that, same to (\ref {convergence_rate_1a}), the contraction factor is the same in the two subdomains and at all grid nodes. However, the factor is related to ${\eta _1}_i^n$ and ${\eta _2}_j^n$, and thus it changes at different time steps. Computation of (\ref {Schwarz_iteration_10}) is made in association with the previous four cases. The contraction factor calculated from the numerical solutions and the theoretical one obtained from (\ref {convergence_rate_3a}) are presented in Table \ref {cla_convergence}, from which it is seen that they are close. The difference is attributed to the linearization adopted in the derivation of the theoretical contraction factor. 

\begin{table}[ht]
  \centering
  \caption{\small Theoretical and numerical contraction factors. (I,J)=(40,80), $\Delta t, \Delta x=0.017$. } 
  \label{cla_convergence}
  \begin{tabular}{|c|c|c|c|c|c|c|}
\hline
\multicolumn{1}{|c|}{Time} & \multicolumn{2}{|c|}{0.085} & \multicolumn{2}{|c|}{0.170} & \multicolumn{2}{|c|}{0.255} \\ \hline
Case & $\bar {\rho }_{theo}$ & $\bar {\rho }_{num}$ & $\bar {\rho }_{theo}$ & $\bar {\rho }_{num}$ & $\bar {\rho }_{theo}$ & $\bar {\rho }_{num}$  \\ \hline
4 & 0.689514 & 0.689539 & 0.690238 & 0.690248 & 0.690612 & 0.690608      \\ \hline  
5 & 0.295124 & 0.295112 & 0.299041 & 0.299082 & 0.302137 & 0.302190      \\ \hline   
6 & 0.457115 & 0.473308 & 0.462127 & 0.478952 & 0.466882 & 0.483229      \\ \hline   
7 & 0.702824 & 0.691392 & 0.703828 & 0.692029 & 0.704226 & 0.692274      \\ \hline  
  \end{tabular}
\end{table}

Now consider the computation of (\ref {Schwarz_iteration_10}) but in association with optimal interface condition (\ref {optimal_int}). Similarly to Proposition 4.2, the following is concluded. 

{\textsc{proposition 5.4}} When interface condition (\ref {optimal_int}) is adopted, the contraction factor for computation of the linearized Burgers equations is
\begin{eqnarray}
\label{convergence_rate_4a}
\bar {\rho }=-\frac {(\alpha +1)({\eta _2}_2^n+\epsilon _2)-R_2(2)}{(\alpha +1)R_1(I-1)+({\eta _1}_{I-1}^n-\epsilon _1)}
            \cdot \frac {(\beta -1)({\eta _1}_{I-1}^n-\epsilon _1)-R_1(I-1)}{(\beta -1)R_2(2)+({\eta _2}_2^n+\epsilon _2))} 
\end{eqnarray} 
\noindent where $I, J\ge 3$. 

{\textsc{Theorem 5.1}} The contraction factor (\ref {convergence_rate_4a}) becomes zero when 
%\begin{eqnarray}
\begin {equation}
\label{optimal_3}
\alpha =\dfrac {R_2(2)}{{\eta _2}_2^n +\epsilon _2}-1, \ \
\beta  =\dfrac {R_1(I-1)}{{\eta _1}_{I-1}^n -\epsilon _1}+1 
\end {equation} 
%\end{eqnarray} 
\noindent in association with $R_1(I-1)R_2(2)+({\eta _1}_{I-1}^n -\epsilon _1 )({\eta _2}_2^n +\epsilon _2)\ne 0 $. In addition, above values of ($\alpha $, $\beta $) is a saddle point of the contraction factor. 

{\it Proof}: The proof can be made by following the steps for Theorem 4.1. 

\begin{figure}[ht] 
\centering
\subfloat[][]{\label{set1}\includegraphics[width=0.4\linewidth]{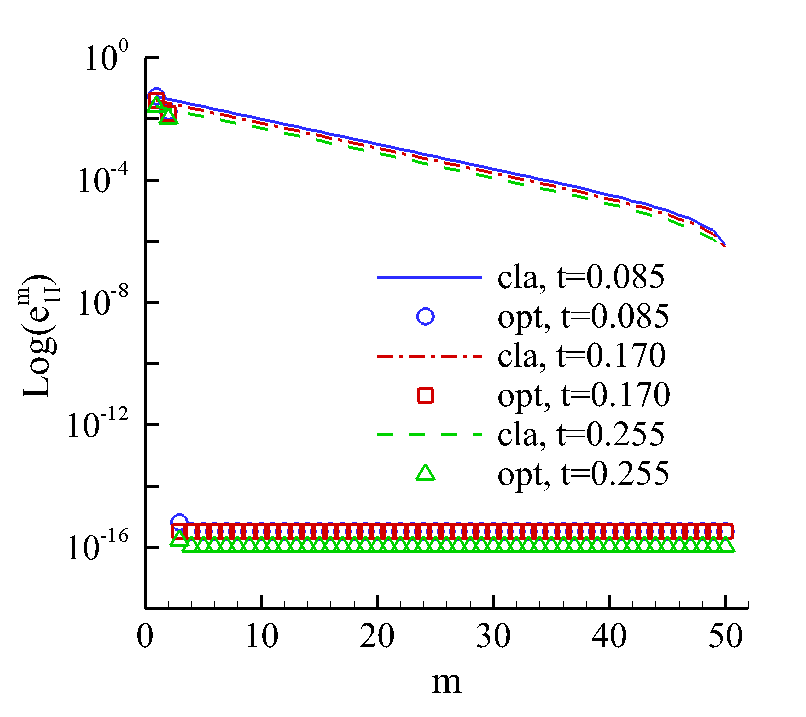}} %\hspace {-0.8 cm} 
\subfloat[][]{\label{set2}\includegraphics[width=0.4\linewidth]{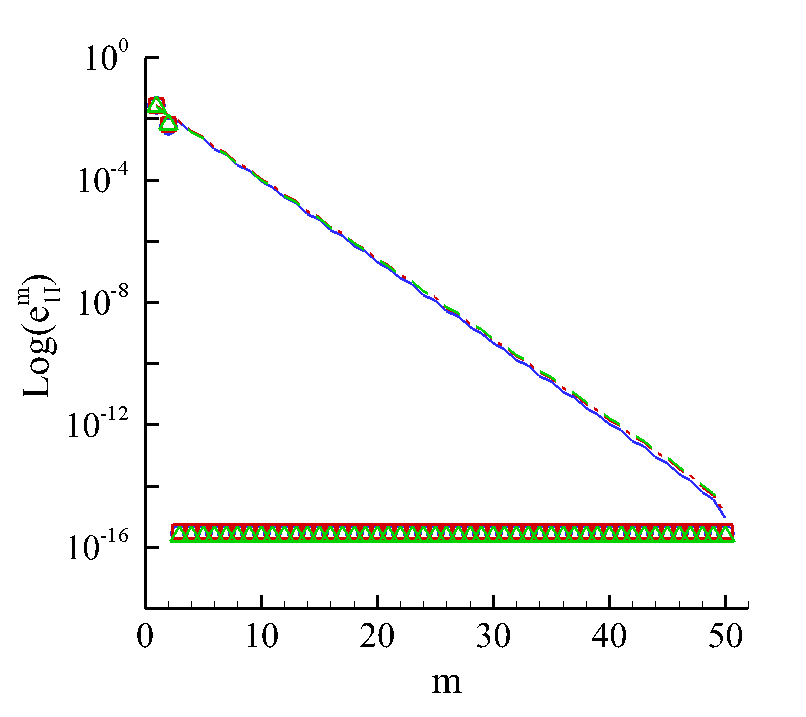}} \\
\subfloat[][]{\label{set3}\includegraphics[width=0.4\linewidth]{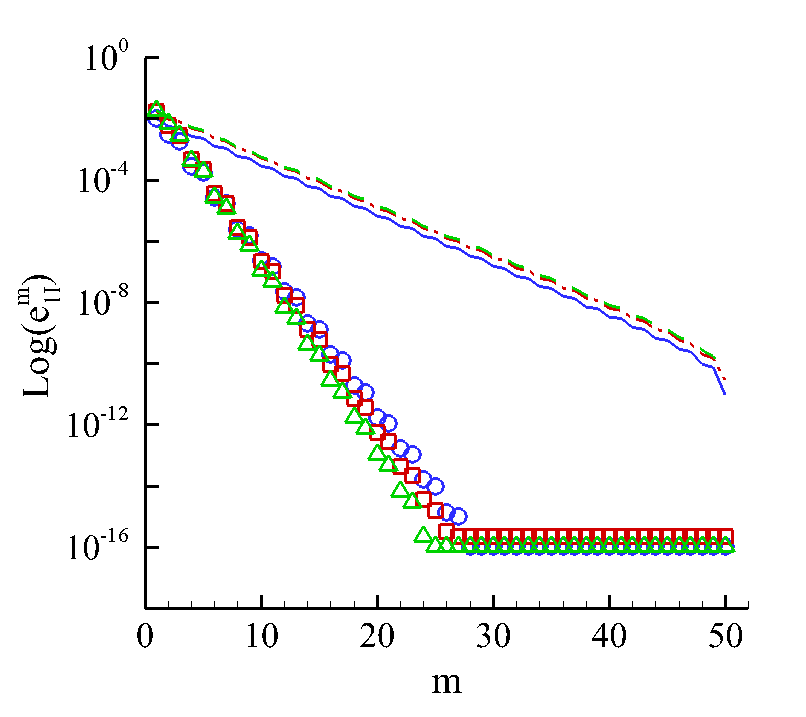}} 
\subfloat[][]{\label{set4}\includegraphics[width=0.4\linewidth]{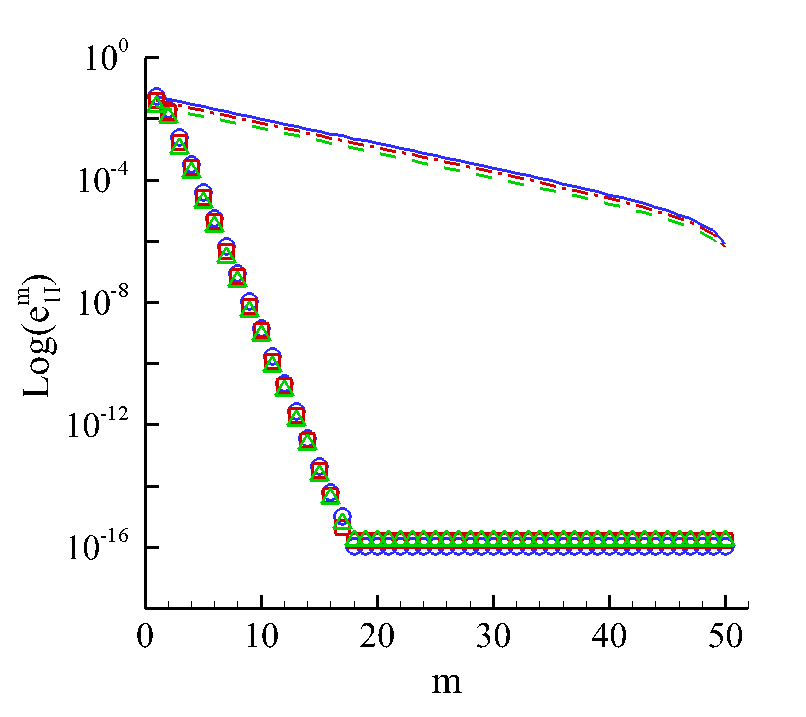}}
\caption{\small Computation with optimal interface conditions for the Burgers equations. a) Case 4. b) Case 5. c) Case 6. d) Case 7.} 
\label{Burgers_opt_convergence}
\vspace {-0.2 cm}
\end{figure}

Numerical results for the four cases in Table \ref {Burgers_cases} are presented in Table \ref {Burgers_opt_convergence_rate}. It is seen that, when the optimal values for ($\alpha ,\beta $) in (\ref {optimal_3}) are adopted, the actual contraction factors in the computation of the first two cases are almost zero, or the theoretical values, while they are somewhat away from zero in the rest two cases. The convergence residual is plotted in Fig. \ref {Burgers_opt_convergence}. The figure shows that, again, computation of the first two cases exhibits the ``perfect convergence'', that is, it converges within about two times of iteration. In the rest two cases, the speed of convergence slows down, and the computation with the optimal interface condition needs many more times of iteration. However, it is still much faster than that associated with the classic interface condition. In the numerical experiment, it is noticed that the slowdown in convergence, as in the last two cases, happens once the diffusion coefficients $b_1$ and $b_2$ become different, even slightly. The reason for the slow down needs further investigation. 

\begin{table}[ht]
  \centering
  \caption{\small Convergence speed in computation of the Burgers equations, together with the optimized interface condition (\ref {optimal_int}). (I, J)=(40,80), $\Delta t, \Delta x=0.017$.} 
  \label{Burgers_opt_convergence_rate}
  \begin{tabular}{|c|c|c|c|c|c|c|c|}
\hline
      Case & Time   & $\alpha$, \qquad \qquad \qquad $\beta $  & $\bar {\rho }_{num}$ \\ \hline
\multirow{3}{*}{4} & 0.085  & $1.90740\times 10^{-1}$, \ $-2.17642\times 10^{-1}$   & $-1.25\times 10^{-14}$ \\ \cline{2-4} 
                   & 0.170  & $1.94634\times 10^{-1}$, \ $-2.12571\times 10^{-1}$   & $ 8.63\times 10^{-15}$ \\ \cline{2-4}  
                   & 0.255  & $1.97387\times 10^{-1}$, \ $-2.09199\times 10^{-1}$   & $ 6.46\times 10^{-15}$ \\ \hline 
\multirow{3}{*}{5} & 0.085  & $6.09587\times 10^{-1}$, \ $-0.99418\times 10^{-1}$   & $ 0.00\times 10^{+0}$ \\ \cline{2-4} 
                   & 0.170  & $6.37715\times 10^{-1}$, \ $-0.37090\times 10^{-1}$   & $-1.69\times 10^{-14}$ \\ \cline{2-4}  
                   & 0.255  & $6.66389\times 10^{-1}$, \ $-9.82950\times 10^{-1}$   & $ 0.00\times 10^{+0}$ \\ \hline 
\multirow{3}{*}{6} & 0.085  & $7.58689\times 10^{-1}$, \ $-4.29258\times 10^{-1}$   & $ 9.44\times 10^{-2}$ \\ \cline{2-4} 
                   & 0.170  & $7.62619\times 10^{-1}$, \ $-1.64360\times 10^{-1}$   & $ 7.72\times 10^{-2}$ \\ \cline{2-4}  
                   & 0.255  & $7.67611\times 10^{-1}$, \ $-9.94238\times 10^{-1}$   & $ 6.31\times 10^{-2}$ \\ \hline 
\multirow{3}{*}{7} & 0.085  & $2.39686\times 10^{-1}$, \ $-1.47483\times 10^{-1}$   & $-1.64\times 10^{-2}$ \\ \cline{2-4}  
                   & 0.170  & $2.44286\times 10^{-1}$, \ $-1.41768\times 10^{-1}$   & $-1.67\times 10^{-2}$ \\ \cline{2-4} 
                   & 0.255  & $2.46275\times 10^{-1}$, \ $-1.39353\times 10^{-1}$   & $-1.68\times 10^{-2}$ \\ \hline 
  \end{tabular}
\end{table}

\vspace{0.3cm} 

\noindent {\bf 6. Concluding remarks}

\vspace{0.2cm} 

A study is presented on the computation of two coupled advection-diffusion-reaction equations. The research starts with the computation of the equations by explicit schemes. Then, it proceeds to an implicit scheme and speedup of convergence by an optimal interface condition. Further, the study is extended to coupled Burgers equations. This study achieves the following main results. 

1. Conditions are presented for convergence of the computation by explicit schemes, and they are easy to use and check in practical problems. Also, such conditions bear clues to speed up convergence. 

2. When an implicit scheme is adopted, an expression for convergence speed is derived. Additionally, an optimal interface condition is presented, which leads to ``perfect convergence", that is, convergence after two times of iteration. 

3. Conclusions drawn in the above linear situations remain mostly true in the situations of the Burgers equations, such as the speedup by the optimal interface condition.

This study sheds light in understanding the computation of coupled advection-diffusion-reaction equations. For instance, this study confirms that convergence becomes slower as grid spacing gets fine, e.g., Fig. \ref {Jacobi_Pseudo_Radius}, which is concluded in computation by waveform relaxation methods. However, it observes that the convergence speed may not increase or even decrease with grid spacing in a wide range of parameters. Additionally, numerical evidence indicates that the overall convergence speed becomes the fastest if an outer iteration is made after every inner iteration, e.g., Table \ref {inner_outer_iteration} and \ref {Burgers_inner_outer_iteration}. Moreover, unlike those obtained in the transformed spaces in the continuous and semi-discretization levels, the convergence speeds' expressions are explicit, and they directly predict those in actual computation. Actually, the method to derive these convergence speeds is extendable to other situations; we have derived the speeds for Poisson equations, which will be reported in a separate paper. 

Some topics deserve further study. One of them is to extend the work of this paper to high dimensions in space, it would be more significant if successful. Another topic will be to understand why the convergence rate with an optimal interface condition remains "perfect" in some scenarios of the Burgers equations, while it slows down substantially in some others. Theoretical analysis will be interesting on the relationship between inner and outer iteration, and why and under what conditions the 1:1 strategy leads to the fastest convergence. It will be interesting to compare the approach in this work,i.e., adopting Schwarz iteration while marching between two adjacent time levels, and that of waveform relaxation. All of these are potential topics of future study. 

\vskip 0.3 true cm

\noindent {\bf Acknowledgments.} This work is supported by NSF (DMS $\#$ 1622453, $\#$ 1622459).  

\bibliography{references}

\end{document}